\documentclass[12pt]{nrc1}
\usepackage{graphicx}
\usepackage{wrapfig}
\usepackage{float}

\usepackage[T1]{fontenc}
\usepackage{fncylab}
\usepackage[english]{babel}
\usepackage{epigraph}
\usepackage{times}

\usepackage{warning}
\usepackage{amssymb}
\usepackage{amsmath}
\usepackage{amsthm}
\usepackage{amscdm}

\makeatletter

\providecommand{\LyX}{L\kern-.1667em\lower.25em\hbox{Y}\kern-.125emX\@}

\usepackage{mdbmaths,mdblogic}
\input headers-notheorem.mtx
\input headers-theorems.mtx

\title{Covers of Abelian varieties as analytic Zariski structure}
\author{\!\!{{M}}.Gavrilovich}
\address{\,\,\,\,\,\,\,{\tiny\sf em\!\!\!\!\!ma\!\!\!ail:\,} gavrilovich at g\!\!\!gm\!\!\!\!\!ma\!\!\!ai\!\!il\!\!l.com}
\def\Bbb{\mathbb}

\def\N{{\Bbb N}}

\def\cl{\operatorname{\textsf{cl}}}

\def\ShafC{{\mathfrak S}}
\def\smooth{{\textrm{sm}}}

\begin{document}

\sloppy 

\maketitle

\begin{abstract}

We use tools of mathematical logic to analyse the notion of a path on an complex algebraic variety, 
and are led to formulate a ''rigidity'' property of fundamental groups specific to algebraic varieties,
as well as to define a bona fide topology closely related to etale topology. These appear as criteria for 
$\aleph_1$-categoricity, or rather stability and homogeneity, of the formal countable language we propose to describe 
homotopy classes of paths on a variety, or equivalently, its universal covering space.

Technically, for a variety $A$ defined over a finite field extension of the field $\Q$ of rational numbers,
we introduce a countable language $L(A)$ describing the universal covering space of $\A(\C) $, or,
equivalently, homotopy classes of paths in $\pi_1(\A(\C) )$. Under some assumptions 
on $A$ we show that the universal covering space of $\A(\C) $ is an analytic Zariski 
structure \cite{ZilberBook}, and present an $\Lww(L(A))$-sentence axiomatising the class
containing the structure and that is stable and homogeneous over elementary submodels.
The ''rigidity'' condition on fundamental groups says that projection of of the fundamental group of a variety
is the fundamental group of the projection, up to finite index and under some irreducibility assumptions,
and is used to prove that the projection of an irreducible closed set is closed in the analytic Zariski structure.

In particular, we define an analytic Zariski structure on the universal covering space
of an Abelian variety defined over a finite extension of the field $\Q$ of rational numbers.


\end{abstract}

\tableofcontents

\section{Introduction }

In \S 1 we describe our approach in a non-technical manner; \S 1.1.1 describes our philosophy
behind the author's thesis \cite{gavr-thesis}, the present paper and \cite{micha-note},
and \S 1.1.2 announces our main results but not detailing definitions. A detailed
exposition of our motivation is found in \S 1.2. In \S 2.1 we give the definitions 
and state the results in \S 2.2. The rest of the paper is devoted to the proof.

\subsection{General Framework}

\subsubsection{Our philosophy}

Is the notion of homotopy on a complex algebraic variety an
algebraic notion? That is, can the notion of homotopy be
characterised in a purely algebraic way, without reference to the
complex topology?

We can restrict to 1-dimensional homotopies only: a 1-dimensional
homotopy is a path, so the question is now whether the notion of a
path on a complex algebraic manifold,  up to fixed point homotopy,
can be characterised in a purely algebraic way.

We provide a partial positive answer to the following more precise
question. Assume that one has an abstract notion of a path up to
homotopy, so that one is able to speak about homotopy classes of
paths, their endpoints, liftings along topological coverings, paths
lying in a subvariety. Can this notion be described without recourse
to the complex topology?

Is it true that one can axiomatise this notion in such a way that
any of its realisations  comes from a choice of an embedding of the
underlying field into $\C$, or equivalently, a choice of a locally
compact Archimedean Hausdorff topology on the underlying field (if
its cardinality is $2^{\aleph_0}$)?

Is the resulting formal theory ``good'' from a model-theoretic point
of view?

Model theory allows a rigorous  formulation of the question as the
problem of proving categoricity of a structure related to the
fundamental groupoid, or equivalently the universal covering space,
of a complex algebraic variety.
 Such
categoricity questions are extensively studied in model theory,
specifically by Shelah \cite{ShelahQA, ShelahQB} and  a short list
of conditions sufficient for categoricity of an
$L_{\omega_1,\omega}$-sentence is known (this is the notion of an
excellent theory).
 Our
model-theoretic analysis shows that the  positive answer to our
question is plausible and is {essentially equivalent to} deep
geometric and arithmetic properties
 of the underlying variety. Some of the properties are known to hold, some
others are conjectured.

We study the interaction between the model theory, arithmetic and
geometry of complex algebraic varieties. Our main results state that
 certain basic  model-theoretic conditions do indeed hold. In general
the proofs {require some technical finiteness and compactness
conditions} and assume some {complex-analytic and arithmetic}
properties and conjectures.
For some classes of varieties, for example Abelian varieties, 
these conditions are known to hold, and 
for these classes
 the results are unconditional. In particular we prove that there exists
an $\aleph_1$-categorical $L_{\omega_1,\omega}$-axiomatisation of
universal covering spaces in such classes.

In \cite[Ch.V]{gavr-thesis}(cf.~also \cite{micha-note})  
we consider a special case where the underlying variety is
an elliptic curve, and prove that the {\bf
natural}  $L_{\omega_1,\omega}$-axiomatisation of the universal
cover of an elliptic curve is $\aleph_1$-categorical; 
analysis there shows that $\aleph_1$-categoricity of that axiomatisation
is essentially equivalent to a arithmetic conjecture 
on Galois representations known for elliptic curves.

Finally we would like to note that the model-theoretic analysis of
universal covers  falls very naturally into the framework of
(analytic) Zariski geometries  started by Hrushovski-Zilber in
\cite{HruZil} and further developed by Zilber and his collaborators
\cite{ZilberZS, ZilberLP,lucy-thesis,martin-zilber-paper,ZilberBook} 
{around an expectation that many basic
mathematical structures may be considered as a model-theoretic
structure with nice properties, above all categoricity}.
Importantly, it has been understood that the model theory relevant
here is essentially non first-order. In fact, our main result is that 
the structures we consider are indeed analytic Zariski as defined in \cite{ZilberBook}, 
thus
providing a series of examples of analytic Zariski geometries.

\subsubsection{Technical summary of  results}
\def\be{\begin{itemize}\item[$\star$]}
\def\ee{\end{itemize}}
In \S\ref{III.3} we define a natural formal
countable language $\Ll$ associated with the universal covering
space $p:\U\ra \A(\C) $ of a complex projective algebraic variety
$\A(\C) $ defined over $\Q$ or $\barq$. 
Assuming subgroup separability of the fundamental group along with its Cartesian
powers, 
we prove that 
\be{\em the positively type-definable sets in $\La$
form a topology analogous to Zariski topology on the set of 
geometric points of a variety,}\ee and, moreover, that

\be{\em the universal covering space $\U^\Ll$, as an $\La$-structure, is 
an analytic Zariski structure \cite[Def.6.1.11]{ZilberBook}}
\ee
By virtue of $\U^\Ll$  being analytic Zariski, we then know

\be{\em the structure $\U^\Ll$ is homogeneous over countable submodels}
($\omega$-model homogeneity), and \emph{realises countably many
types over a countable submodel}.\ee

We then consider in \S\ref{III.5} a fragment of the $\lww(\la)$-theory
$\TwwU$ of $\U^\Ll$ and introduce a natural set of axioms ${\frak X}$
of geometric, analytic Zariski flavour
to show that 

\be{\em the class of models defined by  ${\frak X}$ is stable (in a
non-elementary context) over countable models, and all its models
are homogeneous over submodels.}\ee

These are prerequisites, by Shelah's theory, of categoricity in
uncountable cardinals. Notice that some of the properties, e.g.
atomicity of every model, could, by Shelah's theory, be obtained
just by an $\Lww$-definable expansion of the theory ${\frak X}.$
This, by Shelah's theory, is enough to imply $\aleph_1$-categoricity
of an $L_{\omega_1,\omega}$-class $\Phi$ containing $\U^\Ll$, for an
arbitrary smooth projective variety $\A$ 
with certain conditions on the
fundamental group. (Cf.~Definition~\ref{def:shafarevich} for the
exact definition of the class of algebraic varieties).

Finally we remark that our approach is essentially different
from Zilber's of  \cite{ZilberMTA} since our language $\La$ is in
general stronger than Zilber's. In fact $\La$ ``adjusts'' itself to the
geometric properties of the covering of $A$, and is defined for any
$A$ whereas \cite{ZilberMTA} is restricted to the class of Abelian
varieties.  Our language allows us to produce a sentence in all
cases, {conjecturally categorical for suitably ``self-sufficient''
$A$ whereas \cite{ZilberMTA} is restricted only to considering
Abelian varieties, and those are sometimes obviously not
``self-sufficient'', say Abelian varieties of dimension greater than
1.} We refer to \cite[IV\S6]{gavr-thesis} for details. 
Here we just remark that it is possible to 
consider the language $\La$
corresponding to an ample homogeneous $\Cm$-bundle $A=L^*$ over $X$,
and show that $\La$ defines the 1st Chern class of $X(\C)$ as an
element $c_1\in H^2(\pi_1(X(\C),0),\Z)$ or, equivalently, as an
alternating bilinear Riemannian form $\Lambda\times \Lambda\to {\Z}$.

\subsection{Motivations and implications}

In this section we discuss the motivations behind our choice of the
language and explain our approach in greater detail. In our opinion
the motivations here are more important than the proofs that follow.

We should add that we do not mention yet another motivation relating to category theory 
and Poincare groupoids (\cite[\S I.2.3]{gavr-thesis}, cf.~also \cite{micha-note}), 
as it has no relation to the methods of this paper.

\subsubsection{The Logic approach: What is an appropriate language to talk about paths?
\label{categoriricity:homotopical}}

Abstract algebraic geometry provides a language appropriate to talk
about complex algebraic varieties; what language would be
appropriate to talk about the homotopies on the algebraic varieties,
in particular about paths, i.e.~1-dimensional homotopies? What is
the right mathematical measure to judge appropriateness of the
language for such a notion? 

Abstract algebraic geometry over a field has no complete analogue of
the notion. However, there is a strong intuition based on the naive
notion of a path in complex topology; it is a well-known phenomenon
that naive arguments based on the notion of a path quite often lead
to statements which generalise, in one way or another, to, say,
arbitrary schemes, but which are quite difficult to prove. There
have been many attempts to develop substitute notions, starting from
Grothendieck [SGA1,SGA2,,SGA4$\frac 12$] 
 who developed for this
purpose the notion of a finite \emph{covering} in the category of
arbitrary schemes (\etale{} morphism); see  Grothendieck
\cite{grothendieck-quillen} for an attempt to provide an algebraic
formalism to express homotopy properties of topological spaces, and
Voevodsky-Kapranov~\cite{voevodsky-kapranov} for exact definitions.

Thus, from the point of view of philosophy of mathematics, it is
natural to try to understand why  the notion of a path is so
fruitful and applicable, despite the fact that all attempts to
generalise it to non-topological contexts have had only partial
success.


We intend to propose in this work a model-theoretic structure which
contains an abstract substitute for the notion of a path. The
substitute must possess  the familiar properties of paths appearing
in the topological context, rich enough
 to imply a useful theory of paths; in
particular \emph{they must determine the notion of a path on an
abstract algebraic variety uniquely up to isomorphism}.

 Note that Grothendieck
\cite{grothendieck-quillen}, cf.~also
Voevodsky-Kapranov~\cite{voevodsky-kapranov}, provides a natural
algebraic setup to talk about  paths thereby rather directly leading
to a choice of a language (of 2-functors). Our approach  is in fact
based on a similar idea.

Model theory provides a framework to formulate the uniqueness
property in a mathematically rigorous fashion. Following
\cite{ZilberCovers,ZilberLP} we use the notion of \emph{categoricity
in uncountable cardinals} (of non-elementary classes).  In his
philosophy categoricity is a model-theoretic criterion for
determining when an algebraic formalisation of an object, of perhaps
geometric character, is canonical and reflects  the properties of
the object in a complete way.

In this work we introduce a language $\Ll$ which is appropriate for
describing  the basic homotopy properties of algebraic varieties in
their complex topology, and prove some partial results towards
categoricity and stability of associated structures in that
language. The expressive power of $\Ll$ is studied elsewhere;
here we make the following remarks whose justification can be found 
in \cite[Ch.II]{gavr-thesis}. The language $\Ll$ is capable of expressing properties of
1-dimensional homotopies, i.e.~the properties of paths up to
homotopies fixing the ends. We can speak in  $\Ll$ in terms of
lifting paths to a topological covering, paths lying in closed
algebraic subvarieties (i.e.~a homotopy class has a representative
which lies in the subvariety), paths in direct products and so on.
These properties are sufficient to carry out many basic
1-dimensional homotopy theory constructions. Most notably, following
a construction in Mumford~\cite{MumfordAV} one can definably
construct a bilinear form
$\phi_L:\pii(\A(\C) ,0)\times\pii(\A(\C) ,0)\ra \pii(\Cm,1)$ in the
second homology group $H^2(\A(\C) ,\Z)\cong \bigwedge^2 H^1(\A(\C) ,\Z)$
associated to an algebraic $\Cm$-bundle $L$ over a complex Abelian
variety $X(\C)$. Thus, generally the language has more expressive
power than the one considered originally  by Zilber in
\cite{ZilberMTA}; in particular, some Abelian varieties which are
not  categorical in Zilber's language of \cite{ZilberMTA} are
expected to be categorical in our language. It would be interesting
to know whether our language can interpret Hodge decomposition on
cohomology groups, using the isomorphism $H^n(\A(\C) ,\C)\cong
\bigwedge^n H^1(\A(\C) ,\C)= \bigwedge^n \Hom(\pii(\A(\C) ,0),\C)$
(cf.~\cite{MumfordAV}).
%

The results which we prove towards categoricity in uncountable
cardinalities are partial. We prove categoricity in cardinality
$\aleph_1$ for some special classes of algebraic varieties, e.g. for
elliptic curves.
 We also prove important necessary
conditions, such as stability and homogeneity over models, for much
wider classes.

\subsubsection{The Geometric approach: Analytic Zariski
structures\label{categoricity:analytic}}

The universal covering of an algebraic variety is one of the
simplest analytic structures associated to an algebraic variety and
which is more than an algebraic variety itself; the universal
covering space inherits all the local structure the base space
possesses; and in particular, for a complex algebraic variety it is
a complex analytic space. Thus it is natural to consider it in the
context of Zariski geometries \cite{ZilberZS}: one wants to define a
Zariski-type topology on the universal covering space $\U$ of variety
$\A(\C) $ reflecting the connection between $\U$ and $A$, and such that
$\U$ possesses homogeneity, stability and categoricity properties,
perhaps in a non-first order, $\lww$, way, in a countable language
related to the chosen topology on $\U$.

For this, consider the universal covering space $p:\U\ra \A(\C) $ of an
algebraic variety $\A$. It is natural to assume that the covering map
$p$ and the full algebraic variety structure on $\A(\C) $ are
definable. Then the analytic subsets of $\U$ which are the preimages
$\pinv(Z(\C))$ of algebraic subvarieties $Z$ of $\A(\C) $, are
definable. It is natural to let the analytic irreducible components
of such sets also be definable; one justification for this might be
the desire for an irreducible decomposition.

The above considerations lead us to define a topology on $\U$ as 
generated by 
unions  of analytic irreducible
components of the preimage of  a closed algebraic
subvariety of $\A(\C) $.

It turns out that this topology is rather nice in that it (almost)
admits quantifier elimination down to the level of closed sets, has
DCC (the descending chain condition) for \emph{irreducible} sets,
and can be defined in a countable language {(assuming that the
the Cartesian powers of the fundamental group are subgroup separable, 
a condition we believe 
to be technical
).}  These properties of the topology are
axiomatised in the notion of an analytic Zariski structure in \cite{ZilberBook},
and are 
sufficient to imply the model homogeneity of the structure $p:\U\ra
\A(\C) $, and, more generally, to construct an $\Lww$-class containing
$p:\U\ra \A(\C) $ which is stable over models and whose models are
model homogeneous.  It also turns out that the language obtained in
this way
 is the language appropriate for
describing the paths, as explained in subsection above. We explain
the connection in \S\ref{n:equiv}.

\section{A Zariski topology on a universal covering space of an Abelian variety}


\subsection{Definitions and background}

\subsubsection{Notations and some background}

We briefly introduce basic notions of topology we require. Consult 
\cite[Ch.4,\S\S 2-4]{topology} or \cite{massey}  for details. 


For a 
Hausdorff, locally connected and locally linearly
connected 
topological space $B$ with a distinguished base-point $b\in B$, 
the {\em universal covering space} $(U,u_0)$ of $(B,b)$ is the space of all {\em paths} starting at 
the base-point $b$, i.e. 
continuous maps  $\gamma:[0,1]\longrightarrow  B$, $\gamma(0)=b$, considered up to homotopy
fixing the end-points, and endowed with the natural topology, and further 
equipped with the {\em covering homeomorphism} $p:\U \longrightarrow B$, $p(\gamma)=\gamma(1)$. Two paths 
$\gamma_1,\gamma_2:[0,1]\longrightarrow  B$  are {\em end-point} {\em homotopic}  iff there exists
a {\em homotopy} $\Gamma:[0,1]\times [0,1] \longrightarrow  B$ such that $\Gamma(0,t)=\gamma_0(t),
\Gamma(1,t)=\gamma_1(t), \Gamma(t,0)=\gamma_0(0)=\gamma_1(0), 
\Gamma(t,1)=\gamma_0(1)=\gamma_1(1)$. The {\em fundamental groupoid} $\pi_1(B)$ is the set 
of all paths considered up to end-point homotopy, equipped with the partial operation
of concatenation. A {\em concatenation} $\gamma_0\gamma_1$, $\gamma_0(1)=\gamma_1(1)$ 
of paths is a path which 
first follows the first path $\gamma_0$, and then goes along the second path $\gamma_1$;
this defines concatenation  up to homotopy. The {\em fundamental group} $\pi_1(B,b)=p\inv(b)$ 
is the group of all loops $\gamma\in \pi(B)$, $\gamma(0)=\gamma(1)=b$. A {\em deck transformation}
of $U$ is a homeomorphism $\tau:\U \longrightarrow U$ commuting with $p$, $\tau\circ p = p $. Deck
transformation of $U$ form a group $\Gamma=\pi(U)$ called {\em the deck transformation group}. 
The deck transformation group $\pi(U)$ is canonically identified
with the fundamental group  $\pi_1(B,b)$: to an element $\gamma\in \pi(V')$ 
there corresponds path $p(\gamma_{x'_0,\gamma x'_0})$
where $b'\in U $ is arbitrary such that $p(b')=b$.
 The covering map $p:\U \longrightarrow B$ is a local homeomorphism; a analytic space structure on $B$
induces a unique analytic space structure on $U$.  
 There is a {\em Galois} correspondence between normal subgroups $H<\Gamma$ of $\Gamma$ 
and covering spaces $B^H=^{\rm\!\!\!\!\!\!\! def}U/H \longrightarrow  B$. The map $U\longrightarrow B^H$ is a
universal covering map, and its deck transformation group is $H$; the map
$B^H \longrightarrow B$ is a covering and its deck transformation group is the factorgroup $\Gamma/H$.  

 A map $p:X\ra Y$ is called a \emph{fibration} iff
for any space $Z$ any {\em homotopy} $F:Z\times I\ra Y$ {\em covered} at the
initial time $t=a$, can be {\em covered} at all times $a\leq t\leq b$, $I=[a,b]$ by
some homotopy $G:Z\times I\ra X$ so that $p\circ G(z,t)=F(z,t),
G(z,a)=g(z)$. That is, if map $f(z)=F(z,a):Z\ra Y$ is covered by a
map $g:Z\ra X$, $f(z)=F(z,a)=p\circ g(z),z\in Z$, then there exist a
homotopy $G:Z\times I\ra X$ covering $F:Z\times I\ra X$,
$$ G(z,a)=g(z)$$
 $$F(z,t)=p\circ G(z,t).$$
Homotopy $G$ is called a \emph{covering homotopy with initial
condition $g$}. We also say that homotopy $F$ \emph{lifts} to
homotopy $G$, and that fibration $p:X\ra Y$ has \emph{lifting
property}. We will use extensively the case when $Z=I$ is an 
interval and $g=g(0,a)$ is a simply a point; this case is called
the {\em path-lifting} property. 
A covering is a fibration with discrete fibres.

Often one modifies the definition by restricting $Z$ to a
subclass of spaces, e.g. $Z=I^n$ is required to be a direct product of 
intervals (Serre fibration).  This distinction is not important in this paper.

\subsubsection{Our Assumptions. }

The most interesting, and the only unconditional, example where 
our theorems apply,  is that of  $\U/\Gamma$  an Abelian variety:
$\U=\C^{2g}$, $\Gamma=\pi(\U/\Gamma,0)$ is a lattice in $\U$.  

However, our assumptions are geometric; in particular the assumptions do not mention the 
group structure of an Abelian variety. We call the corresponding 
class of varieties {\em LERF}. 

We assume  $\U$ is a smooth 
complex analytic space equipped with an 
free cocompact action $\Gamma : \U \longrightarrow  \U$ of a 
subgroup separable (cf.~\ref{sec:lerf})
finitely generated group $\Gamma:\U\longrightarrow \U$. Further 
we assume that all Cartesian powers of $\Gamma$ are 
subgroup separable, and that $\U/\Gamma$ is a projective algebraic
variety.

\paragraph*{ {\label{sec:lerf}Subgroup separability of $\pi(\U)$}.}
A group $\Gamma$ is called {\em subgroup separable}, or {\em locally
extended residually finite}, often abbreviated {\em lerf}, iff for
any finitely generated subgroup $G<\Gamma$ and an element $g\not\in G$ there exists a {\em
finite} index subgroup $H$ such that $G<H$ and $g\not\in H$. This is
a non-trivial property rather hard to establish; it is known that
the fundamental groups of complex curves
(\cite{surface-groups-are-geometric}) and $\Z^n$, $\SL_2(\Z)$ are
subgroup separable; however, it is known that $F_2\times F_2$
(\cite{miller}) is not subgroup separable, and so in general the
products of subgroup separable groups are not subgroup separable.
This property may be reformulated topologically: the group
$\Gamma=\pi_1(A)$ is subgroup separable if and only if for any finitely
generated $G < \Gamma$ and any compact subset $C\subset A^G=U/G$, 
the covering splits as $A^G \longrightarrow  A^H \longrightarrow  A$ such that $A^H \longrightarrow  A$ is a 
finite covering and the compact $C$ maps to $A^H$ by a homeomorphism.
In fact, we need this
property only when $G$  is  the fundamental group of an algebraic
subset of $\A$. 

\paragraph*{ LERF varieties.}

The above enables us to define the class of LERF varieties to which 
our theorems apply. 

\begin{defn}\label{def:shafarevich}
We call a smooth projective algebraic variety $\\A(\C) $ \emph{LERF}
if
all finite Cartesian powers of the group of deck transformations 
$\pi(\U)$ are subgroup separable.  
\end{defn}


\subsubsection{{Co-etale topology, its core and inner core} }

We define topologies on $\U$ and Cartesian powers on $\U$.

\paragraph*{ Definition of the co-etale topology.}
We give 3 equivalent definitions of {\em co-etale topology} on U; 
we  prove the  equivalence in Decomposition Lemma~\ref{lem:noeth:ana}.

\begin{defn}\label{def:Tt} 
(I) A subset of $\U^n$, $n>0$, is {\em closed in co-etale topology} $\ShafC$ iff 
it is either 
	(I)(i) an irreducible analytic component of a closed analytic set 
	such that the set is set-wise 
        invariant under the action of the fundamental group,
	or
	(I)(ii) a closed analytic set such that each of its analytic irreducible component 
	satisfies (I)(i) above.

We call a closed analytic subset $Z$ of $\U^n$ {\em unfurled} iff every connected component of $\U$ is irreducible. 
It is known  that every smooth closed analytic set is unfurled. 

$(C)$ A subset of $\U^n$, $n>0$, is {\em closed in co-etale topology} $\ShafC$ iff
it is either
        $(C)(i$) a connected component of an unfurled closed analytic set 
	such that the set is 
	set-wise invariant 
		under action of a finite index subgroup of the fundamental group,
	or
        $(C)(ii)$ a closed analytic set such that each of its analytic irreducible components  
        satisfies (C)(i) above.
 
$(C')$  A subset of $\U^n$, $n>0$, is {\em closed in co-etale topology} $\ShafC$ iff
it is either
        $(C')(i)$ a connected component of an unfurled closed analytic set such that the set 
	is set-wise invariant  
                under action of a finite index subgroup of the fundamental group,
	or
        $(C')(ii)$ a countable intersection of sets  as in $(C')(i)$ 
\end{defn}

\paragraph*{ Countable core $\mathsf C_0$.}
By our assumptions, $\A(\C) =\U/\Gamma$ is a complete projective algebraic variety defined over $\barq$, 
and therefore by Chow's Lemma every closed analytic subset of $\A(\C) $ is in fact algebraic 
and defined over a finitely generated subfield of $\C$. This enables us to speak of 
{\em the field of definition of a $\Gamma$-invariant closed analytic subset} of $\U/\Gamma$, 
as $\Gamma$-invariant closed analytic subsets are in 1-1 correspondence with closed analytic subsets of $\A(\C) $.

This enables us to define the following. 

\begin{defn}  The {\em countable core}  $\mathsf C_0$ consists of closed sets that are unions of 
irreducible components of $\Gamma$-invariant closed sets defined over $\barq$.  \end{defn}Note that 
a point $u\in \U$  is in the countable core iff $p(x)\in A(\barq)$. 

In Lemma~\ref{sec:la-def:lem} we prove that core sets are enough to define all sets; 
in the following  way: 
that every irreducible co-etale closed subset $Z\subset \U^n$  
can be represented as 
a connected component $Z\times\{g\}$ of a hyperplane section 
$ Z' \cap U^n\times \{g\}$ of a co-etale closed set $Z'$ in the countable core. 
 

\paragraph*{ {Countable inner core $\mathsf C_\emptyset$} .}


In fact, in our structure we may define analogs of sets over $\Q$ (or 
perhaps the maximal Abelian extension of $\Q$), and not just $\barq$.

\begin{defn}  The {\em countable inner core} $\mathsf C_\emptyset$ consists of 
the subsets of $\U^n\times U^n$ defined by 
relations $x'\simH y'$ and $x'\simz y'$ where 
$Z\subset \A^n$ is a closed subvariety defined over the field of definition of $\A$, $H$ 
a finite index subgroup of $\Gamma$, 
and the relation is defined as follows. 
\bi\item[] $x'\simz y'$ $\iff$ points $x'\in\uU^n$ and $y'\in\uU^n$
lie in the same (analytic) irreducible component of the
$\Gamma$-invariant closed analytic set $\pinv(Z(\C))\subset\uU^n$.
\item[] $x'\simH y'$ $\iff$ $\exists \tau\in H^n :\tau x'=y'$.
\ei
We shall also consider 
\bi \item[] $x' \sim_{Z,\etAh}^c y'$ iff $x'$ and $y'$ lie in the same
 \emph{connected} component of the preimage $\phh\inv (Z_i(\C))$,
 $Z_i\subset \etAh(\C)^n$  an irreducible component of algebraic
 variety $\pHh\inv(Z(\C))\subset \etAh(\C)^n$.\ei \end{defn}

\subsection{Our Results: Definition of analytic Zariski structure, and the main theorem. }

We have defined a topology on every Cartesian power of U, and the notion of countable core. 

Every co-etale closed set is closed in analytic topology, and thus possesses
the dimension; let this be the dimension function of the analytic Zariski structure.

\begin{thm} The data as defined above, form an analytic Zariski structure as 
defined in \cite[Def.6.1.11]{ZilberBook}. Moreover, the analytic Zariski structure belongs to
an explicitly axiomatised $\Lww$-class $\Ax(\A(\C) )$ that is $\omega$-stable over submodels, every 
model is $\omega$-homogeneous. \end{thm}
\begin{cor}[] Every countable model extends uniquely to a model 
of cardinality $\aleph_1$.  It is consistent with ZFC that every countable model
extends uniquely to a model of cardinality continuum. \end{cor}


The rest of paper is devoted to the proof of these claims; see \S6, Theorem~\ref{prp:main} and 
Theorem~\ref{Th:1}.

We also formulate a conjecture; see \cite[\S IV.6-\S IV.7]{gavr-thesis} or 
a forthcoming paper to clarify 
its relationship to a categoricity conjecture of \cite{ZilberMTA}.

\begin{conj} For generic complex Abelian varieties $\A$ defined over a number field, 
an analogously defined $\lww$-class $\Ax(\A(\C) \times \Cm)$  
 is analytic Zariski, excellent 
and categorical in uncountable cardinalities. A sufficient condition
is that the Mumford-Tate group of $\A$ is the symplectic group, i.e.
the largest possible.
\end{conj}






\subsection{Reduction to unfurled subsets: equivalence of the  definitions 
\label{sec:noeth}}

In this section we prove that the definitions \ref{def:Tt}$(I)$ and \ref{def:Tt}$(C)$  of the collection $\ShafC$ do agree. It is the main prerequisite to prove that $\ShafC$ is a topology. 

\subsubsection{Prerequisites on analytic irreducible decomposition and coverings in algebraic geometry}

\paragraph*{ Irreducible Decomposition in smooth analytic spaces.}

To avoid confusion, below we say ``an open ball'' to mean a
neighbourhood open in complex topology, not in the analytic Zariski topology.

\begin{fact}\label{fact:local} Let $\U$ be a 
smooth complex analytic space, and let
$Y,Z\subset \U$ be closed analytic subsets in $\U$.
 Then \bi

 \item (irreducible decomposition) $Z$ admits a unique decomposition $Z=\cup_{i\in\N}Z_i$ 
into a countable union of analytic irreducible closed subsets $Z_i$'s. 

 \item (analyticity is a local property) a set $X\subset \U$ is analytic
 iff for all $x\in X$, there exists an open ball $x\in B_x$
 such that $X\cap B_x$ is an analytic subset of $B_x$
\item (local identity principle) for an open ball $B\subset \U$, if $Y$ is irreducible and
$Y\cap B\subset Z\cap B$ then $Y\subset Z$

\item (local identity principle; analytic continuation)
for an open ball $B\subset \U$, if $Y$ and $Z$ are
irreducible, and $Y\cap B$ and $Z\cap B$ have a common irreducible
component, then $Y=Z$

\item (density of smooth points)  for an open ball $B\subset \U$,
if $Z_0\subset Z\cap B$ is an irreducible component of $Z\cap B$,
then there exist a point $z_0\in Z_0$ and an open ball $z_0\in
B_0\subset B$ such that $B_0\cap Z\subset Z_0 $

\item  (local finiteness) a compact set
$C\subset \U$ intersects only finitely many irreducible components
of a closed analytic set $Z$

\item (analyticity of a union of irreducible components)\label{union:analytic}
 a union of, possibly infinitely many, irreducible components of
an analytic set is analytic

\item (irreducible decomposition)
if $Y\subset Z$ and $Y$ is irreducible, then $Y$ is contained in
an irreducible component of $Z$

\item (smooth points of irreducible sets) 
the set of smooth points of an irreducible set is connected; 
consequently, the irreducible decomposition $Z=\cup_i Z_i$ of a closed analytic set $Z$ 
is determined by the decomposition $Z^{sm}=\cup_i(Z^{sm}\cap Z_i)$ into connected components 
of the set of its smooth points. 
\ei
\end{fact}

\bp Those are well-known properties of smooth complex analytic spaces.

(1) is by \cite[\S5.4,Theorem, p.49]{cirka}.
By Prop.~5.3 of \cite{cirka}, Theorem 5.1 [ibid.] states (7) and
(6). Corollary 2 of Prop.~5.3 [ibid.] implies (3) and (4). Theorem
5.4 [ibid.] implies (5). (2,3,4) together imply (8).
(9) is by [17,\S5.4,Theorem].
\ep


\paragraph*{ Finite topological coverings in algebraic geometry.}
We also need a form of  Riemann existence theorem. 

\begin{fact}[Generalised Riemann existence theorem]\label{gaga}\label{GAGA} Let $\AC$ be a normal algebraic variety over $\C$.
If $q:T\ra \A(\C) $ is a finite covering of topological spaces, then $T$
admits a structure of a complex algebraic variety such that
$q_\mathrm{top}:T\ra \AC$ becomes an algebraic morphism, i.e. there
exists an algebraic variety $B(\C)$ over $\C$, an algebraic morphism
$q_\alg:B(\C)\ra \AC$, and a homeomorphism $\phi:T\ra B(\C)$ of
topological spaces such that the diagramme of topological spaces
commutes
\becd T @>q_\mathrm{top}>> \AC \\
      @V \phi VV  @V\id VV\\
      B(\C) @>q_\alg>>\AC
\eecd

Moreover, the homeomorphism $\phi:T\ra B(\C)$ is well-defined up to an automorphism of
$B$ commuting with the covering morphism $q_\alg$.
\end{fact}

\bp 
Grothendieck [SGA1,Exp.XII,Th.5.1]; by a variety over $\C$ we mean a Noetherian scheme 
of finite type over $\C$. One may also
look in \cite[Appendix B,\S3,Theorem 3.2]{HartAG} for some
explanations.

\ep

\subsubsection{Reduction to unfurled subsets : the  proof  }

For a subset $Z\subset \U$, let $\Gamma Z
=\bigcup\limits_{\gamma\in\Gamma} \gamma Z'$ denote the $\Gamma$-orbit of
set $Z$.

For $H\nrmfin\Gamma$, let $\phh:\uU\ra \uU/\simH$ be the factorisation
map since $\A=\U/\Gamma$; by Fact~\ref{gaga}, we choose and fix isomorphisms
$\etAh(\C)\cong \factor \uU {\,\simH}$ where $\etAh(\C)$ is an
algebraic variety; the deck group of covering $\etAh(\C)\ra \A(\C) $
is the finite group $\Gamma/H$.

\begin{lem}[First Decomposition lemma; Noetherian property; Reduction to Unfurled Subsets]\label{lem:noeth:ana}
Assume $\A$ is LERF.

Every $\Gamma$-invariant analytic closed set has a decomposition as a finite union
of unfurled closed analytic subsets invariant under the action of a finite index 
subgroup of $\Gamma$.

In other words,
a $\Gamma$-invariant analytic closed set has an analytic decomposition
of the form
$$ W'= H Z'_1\cup \ldots\cup H Z'_k,$$
where $H\nrmfin \Gamma$ is a finite index normal subgroup of $\Gamma$, the
analytic closed sets $Z'_1,\ldots,Z'_k $ are irreducible, and for
any $\tau\in H$ either $\tau Z'_i=Z'_i$ or $\tau Z'_i\cap
Z'_i=\emptyset$. 

Such decomposition also exists for closed analytic sets invariant
under the action of a finite index subgroup of $\Gamma$.

\end{lem}

\bp Let us prove that $(a)$ there exists a decomposition as above
without the condition on intersections, and then prove $(b)$ the
irreducible components satisfy $\tau Z'_i=Z'_i$ or $\tau Z'_i\cap
Z'_i=\emptyset$ for $\tau\in\Gamma$.

The proof of $(a)$ is relatively simple, and follows from the
Fact~\ref{fact:local} in a rather straightforward way; we do it
first.

The proof of the second claim $(b)$ uses rather more delicate local
analysis of the structure, and several local-to-global properties of
analytic subsets of smooth complex analytic spaces as well as some finiteness
properties of Zariski geometry of algebraic varieties.

So let us start to prove $(a)$. Let $Z'$ be an irreducible component
of $\pinv (Z(\C))$; by $\Gamma$-invariance of $\tZ$, for any $\gamma\in
\Gamma$, the set $\gamma Z'$ is also an irreducible component of $\tZ$,
and so $\Gamma Z'$ is a union of irreducible components of $p\inv
Z(\C)$; thus, by Fact~\ref{fact:local} above, $\Gamma Z'\subset \tZ$ is
analytic.

The covering morphism $p:\U\ra \A(\C) $ is a local
isomorphism, and analyticity is a local property; by
$\Gamma$-invariance of $\Gamma Z'$, it implies $p(\Gamma Z')$ is analytic.
For different irreducible components $Z'_1\neq Z'_2$ of $\tZ$ it can
not hold that $p(Z'_1)\subsetneq p(Z'_2)$; indeed,  then $\Gamma
Z'_1=\pinv p(Z'_1)\subset \Gamma Z'_2=\pinv p(Z'_2)$, and so
$Z'_1=\bigcup(Z'_1\cap \gamma Z'_2),\gamma\in \Gamma$; thus, $Z'_1$ can
not be irreducible unless $Z'_1\subset \gamma Z'_2$, for some
$\gamma\in\Gamma$. To conclude, closed sets $p(Z')$, $Z'$ vary among
irreducible components of an algebraic subvariety $Z(\C)$, cover the
whole of $Z(\C)$; they are also irreducible. Thus they are the
analytic irreducible components of $Z$. The
analytic irreducible components of an algebraic set are algebraic
and irreducible  by \cite{Rossi},  and thus they are the algebraic irreducible
components; in particular there are only finitely many of them. That
gives the required decomposition.


Now let us start to prove $(b)$. First of all, note that we
may suppose $Z$ to be irreducible.

\def\Zddva{\bigcup\limits_{Z'_i\neq Z'_j} Z'_i\cap Z'_j
}
\def\Zdnn{\bigcup  Z'_{i_1}\cap\ldots\cap Z'_{i_n}
}
\def\Zdn{Z'^{(n)}}

Let $\Zdn=\Zdnn$ be the union of all intersections of $n$-tuples of
different irreducible components of $\pinv(Z(\C))$.
\begin{claim}\label{claim:216} The set $p(\Zdn)$ is
an algebraic subset of $Z(\C)$, for $n>0$. For $n$ sufficiently
large, $\Zdn$ is empty.
\end{claim}

\bp By the local finiteness (Fact~\ref{fact:local})  a compact
subset intersects only finitely many of the irreducible components
$\gamma Z'_i$'s; thus $\Zdn$ is locally a finite union of
intersections of analytic sets, and therefore is analytic. By the
$\Gamma$-invariance of $\gamma Z'_i$'s it is $\Gamma$-invariant, and thus
$p$ provides a local isomorphism of $\Zdn$ and its image; therefore
the image $p(\Zdn)$ is analytic. By Chow Lemma 
this implies it is in fact algebraic. If $n$ is greater then the number
of local irreducible components at a point of $Z$ in $\A$, then by
Fact~\ref{fact:local}(local identity principle) $\Zdn$ has to be
empty.\ep

The claim above implies $\Zdn$ are co-etale closed, for any $n$. By
Claim $(a)$ of Lemma, we may choose finitely many points $z'_i$'s so
that any irreducible component of $\Zdn$, for each $n>0$, contains a
$\Gamma$-translate of one of $z'_i$'s.

By Fact~\ref{fact:local}(5) every point $z'_i$ is contained in only finitely many
irreducible components of $\p\inv(Z(\C)$. Let  $Z'_1,\ldots,Z'_k$ be 
all the irreducible components of $\p\inv(Z(\C)$ containing at least one 
 of the points $z'_i$'s.

For a subset $V\subset \U^n$, define {\em the deck transformation group of $V$} as 
$\pi(V)= \{\gamma\in\Gamma^n: \gamma V\subset V\}$.  If $V$ is a connected component 
of $\Gamma$-invariant set $\pinv(p(V)$, then $\pi(V)$ is canonically identified 
with the fundamental group $\pi_1(V(\C),x_0)$: 
to an element $\gamma\in \pi(V')$ there corresponds path $p(\gamma_{x'_0,\gamma x'_0})$
where $x'_0\in V'$ is arbitrary such that $p(x'_0)=x_0$.

Notice that $\pi(Z'_i)=\pi(Z'_i\cap (\Gamma Z'_i)^\smooth)$ where $(\Gamma Z'_i)^\smooth$ is
the set of smooth points of $\Gamma Z'_i$, and that by Fact~\ref{fact:local}(9)  
the set $Z'_i\cap (\Gamma Z'_i)^\smooth)$ is a connected component of  $(\Gamma Z'_i)^\smooth$.
By the topological argument above, $\pi(Z'_i)$ is the fundamental group of a 
constructible algebraic set $p(Z'_i)^\smooth$. As a constructible algebraic set,
it admits a finite triangulation into simplices, e.g. by o-minimal cell decomposition, 
and this implies that its fundamental group is finitely presented. In particular,
it is finitely generated and we may apply subgroup separability of $\Gamma$ to find
 a normal finite index subgroup $H\subset \Gamma$ such that
$H Z'_i\neq HZ'_j$ for
$i\neq j$, i.e.~$p_H(Z'_i)\neq p_H(Z'_j)$. 

Consider $Z'_i\cap hZ'_i,h\in H$ and assume 
$\emptyset\subsetneq Z'_i\cap hZ'_i\subsetneq Z'_i$. Then
there exists $\gamma\inv\in\Gamma$ such that $\gamma\inv z_j'\in
Z'_i\cap hZ'_i$, i.e.~$z'_j\in \gamma Z'_i \cap \gamma h Z'_i =
\gamma Z'_i \cap  h'\gamma Z'_i$. Both $\gamma Z'_i$ and $h'\gamma
Z'_i$, $h'\in H$ are connected components containing $z'_j$ and by
definition we have chosen $H$ small enough so that $H \gamma Z'_i
\neq H h' \gamma Z'_i$, a contradiction.

In other words, we have proven that there exists a normal finite index
subgroup $H<\pi(\A(\C) $ such that $Z'_i$ is a connected component of
$\phh\inv\phh(Z'_i)$, i.e.~the connected components of the preimages
of the irreducible components of $\phh p\inv(Z(\C))$ are
irreducible. \ep

\subsubsection{\label{n:equiv} Equivalence of the two definitions of co-etale topology}

The next corollary shows equivalence of the two definitions
of co-etale topology.

Notice that the notion of an $H$-invariant set is essentially
algebraic: an $H$-invariant set is a preimage of a closed algebraic subset in the finite
cover $\etAh(\C)$. Thus, the meaning of the next corollary that in
fact co-etale closed sets encode a mix of algebraic data and
\emph{topological, homotopical} data, not of analytic one.

\begin{cor}\label{cor:union.con}
Definitions~\ref{def:Tt}(I) and ~\ref{def:Tt}(C) are equivalent.
In particular, an irreducible co-etale closed set is a connected component 
of a unfurled closed analytic set invariant under action of a finite index subgroup
of the fundamental group.
\end{cor} \bp
Lemma~\ref{lem:noeth:ana} above implies that each co-etale irreducible closed set according to 
~\ref{def:Tt}(I) is also closed according to ~\ref{def:Tt}(C), i.e.
is a  connected component of a a unfurled closed analytic set invariant under action of a finite index subgroup
of the fundamental group.
.

On the other hand, the lemma implies that each $H$-invariant set is
a finite union of sets of the form $HZ'_i$ where $Z'_i$ are
irreducible. Then, $\Gamma Z'_i$ is also closed analytic as a finite
union of translates of $H Z'_i$, and moreover, each translate of
$Z'_i$ is an irreducible component of $\Gamma Z'_i$ and thus co-etale
closed. This implies every (C)-closed set is also (I)-closed. \ep

\paragraph*{ An algebraic reformulation.}
The Lemma has the following algebraic consequence. All the notions
mentioned in the Corollary are preserved under replacing the ground
field by another algebraically field; 
thus it holds for any characteristic 0 algebraically closed field
instead of $\C$. One may think of this property as a rather weak
property of irreducible decomposition for the \emph{co-etale}
topology; it is also a statement about a resolution of non-normal
singularities.

\begin{cor}\label{cor:decomp.alg} Let $\A$ be LERF. Then for any
closed subvariety $Z\subset \AC$, there exists a finite \'etale
cover $q:\etAhc\ra \AC$ such that, for any further \'etale cover
$q':\etAg(\C)\ra \etAh(\C)$, the connected components of
${q'}\inv(Z_i)\subset \etAg(\C)$ are irreducible, where $Z_i$'s are
the irreducible components of $q\inv(Z)$.
\end{cor}\bp
Indeed, it is enough to take $H$ as in Decomposition Lemma.\ep

Note that when $Z$ is normal, the corollary is a well-known
geometric fact.

\subsection{Co-etale  topology is a topology.  }

\begin{lem} \label{lem:topology}\label{cor:conn.comp}\label{lem:decomp.ana} 
(a) The collection $\ShafC$ of subsets of $\U^n$ forms a topology, for every $n$. 
(a$'$) Moreover, the collection $\ShafC$ satisfies Axioms (L1)-(L8) of \cite{ZilberBook}.
(b) An $\ShafC$-irreducible $\ShafC$-closed set is analytically irreducible closed set. 
(c) An analytically irreducible component of a $\ShafC$-closed set is 
$\ShafC$-closed $\ShafC$-irreducible.  \end{lem}\bp 
(b)  By Definition~\ref{def:Tt}(I), a co-etale
irreducible co-etale closed set $W'$ is a countable union of
irreducible component of $\Gamma$-invariant closed analytic sets. Those
components are co-etale closed by definition, and thus co-etale
irreducibility implies the union is necessarily trivial. Thus, the
set is an analytic irreducible component of a $\Gamma$-invariant set,
i.e.~in particular irreducible as an analytic set.

(c) is immediate by Definition~\ref{def:Tt}(I). 

(a) As  $\ShafC$ consists only of closed analytic sets,  an analytic irreducible component of a finite union of
$\ShafC$-closed sets is an analytic irreducible component of one of them; this shows that $\ShafC$ is closed under finite 
union. To prove $\ShafC$ is closed under infinite intersection, we first observe that an irreducible component 
of an infinite intersection (that is still a closed analytic set) is necessarily the intersection of irreducible
closed analytic components of these sets; by the descending chain condition for analytic irreducible 
closed sets, the intersection is necessarily finite. Thus, by Definition~\ref{def:Tt}(I) it is enough to show
that each irreducible component of the intersection of irreducible $\ShafC$-closed sets is $\ShafC$-closed (irreducible).

Thus, it is enough to prove that the intersection of two 
irreducible $\ShafC$-closed sets, say X and Y, is $\ShafC$-closed. Now, by Definition~\ref{def:Tt}(C), X and Y are 
{\em connected} components of closed analytic sets $X'$ and $Y'$ invariant under action of finite index subgroups, say
H and G, of the fundamental group. Then, $X\cap Y$ is a {\em connected} component of the intersection 
$X'\cap Y'$ that is invariant under action of $H\cap G$, the latter also being a finite index subgroup. 
By Definition~\ref{def:Tt}(C), this implies that $X\cap Y$ is $\ShafC$-closed. This proves (a); note the interplay
between $(I)$ and $(C)$ of Definition~\ref{def:Tt}.  

(a$'$) We have just proven (L1); axioms (L2-L7) are immediate by inspection of any of the definitions.  
(L8) requires ~\ref{def:Tt}(C) : a hyperplane section of a {\em connected} component of a closed analytic set 
invariant under action of a finite index subgroup is a {\em connected} component of the intersection 
that is also invariant under a finite index subgroup. (This argument does not work for irreducible
components, as they may intersect). 

\ep 

\subsection{Good dimension notion : (DP), (DU), (SI), (AF)}

The following properties are defined in \cite[\S3.1]{ZilberBook}. 
Following notation there, $S \subseteq_{cl}  S'$ reads $S$ is a closed subset of $S'$,
$S \subseteq_{an}  S'$ reads $S$ is an analytic subset of $S'$, and $S \subseteq_{op}  S'$
reads $S$ is an open subset of $S'$.

\begin{lem}[Good dimension]

(DP) Dimension of a point is 0 

     (DU) Dimension of unions: $\dim(S_1  \cup  S_2  ) = \max(\dim S_1  , \dim S_2 )$
    
 (SI) Strong irreducibility: For $S \subseteq_{cl}  V \subseteq_{op}  \U^n$ , $dim S_1  < dim S$, if $S$
is irreducible and $S_1  \subseteq_{cl}  S$ is closed, then $S_1  = S$

     (AF) Addition formula: For any irreducible $S \subseteq_{cl}  V \subseteq_{op}  \U^n$ and a
projection map $pr : \U^n \longrightarrow  \U^m$ ,
                $$\dim S = \dim \pr (S) + \min_{ a\in pr (S)} dim(\pr\inv (a) \cap  S).$$
                                    
   (PS) Presmoothness: For any closed irreducible $S_1  , S_2  \subseteq  \U^n$, the
dimension of any irreducible component of $S_1  \cap  S_2 $ is not less than
                       $$ \dim S_1  + \dim S_2  - \dim \U^n .$$

\end{lem}
\bp These are inherited from complex analytic geometry. \ep

\subsection{Analyticity (AS), (SI),(DP),(CU), (INT),(CMP), (CC)}



Recall that \cite[\S6.1.2]{ZilberBook} distinguishes a class of sets in a topology that
he calls 'analytic'. Namely, in a topology $T$ a locally closed set  $S$  
is called {\em analytic} in an open set $U$ iff $S$ is a closed subset of $U$ and for every $a\in S$
there is an open $a\in V_a \subseteq_{op}  U$ such that $S \cap  V_a$ is the union of
finitely many relatively closed irreducible subsets. Note that by Fact~\ref{fact:local}(6,7), 
a locally closed analytic set is analytic in this sense: take $V_a$ to be the 
completement of the union of the irreducible components of $S$ not containing $a$. 
This argument also works for co-etale topology, i.e., in co-etale topology,
each locally closed set is analytic in this sense.

Next Lemma establishes (INT), (CMP),(CC) and (AS) of [loc.cit., \S6.1],
and therefore, that $\U$ is a {\em topological structure with a good dimension theory} 
[loc.cit.,Def.6.1.1].

\begin{lem}[Analytic sets] 
(INT) (Intersections) 
If $S_1  , S_2  \subseteq_{an}  \U^n$ are irreducible and analytic in $\U^n$, then $S_1  \cap  S_2 $ is analytic in $\U^n$ 

(CMP) (Components) If $S \subseteq_{an}  \U^n$ and $a \in  S$ then there is $S_a \subseteq_{an}  \U^n$, a finite
       union of irreducible analytic subsets of $\U^n$, and some $S'_a \subseteq_{an}  \U^n$ such that
       $a \in  S_a \setminus S'_a$ and $S = S_a \cup  S'_a$ 

(CC) (Countability of the number of components) Any $S \subseteq_{an}  \U^n$ is a
     union of at most countably many irreducible components.

   (AS) [Analytic stratification] Every locally closed set is analytic.

   (aPS) [Analytic Presmoothness] If $S_1  , S_2  \subseteq_{an}  V \subseteq_{op}  \U^n$ and both $S_1  , S_2 $ are 
irreducible,
then for any irreducible component $S_0$ of $S_1  \cap  S_2 $
                     $$\dim S_0 \geq  \dim S_1  + \dim S_2  - \dim \U^n.$$

\end{lem}
\bp Immediate by Fact~\ref{fact:local}.
\ep

\subsection{{$\Theta$-definable sets, $\Theta$-generic points and $\Theta$-definable closure}}

Recall that $\uU/\Gamma\cong \AC$ has the structure of an algebraic
variety over $\C$ and that the $\Gamma$-invariant sets are in a bijective correspondence
with the algebraic subvarieties of $\\A(\C) $. This suggests us that
we may try to pull back to $\uU$ the notion of a generic point in $\AC$.

The following definition behaves well only for $\Theta\subset\C$
algebraically closed.

\begin{defn} We say that a $\Gamma$-invariant co-etale closed subset
$W'\subset\uU$ is \emph{defined over an algebraically closed subfield $\Theta\subset\C$}
iff $p(W')\subset\AC$  is a subvariety defined over $\Theta$.

An co-etale  closed set is \emph{defined over a subfield
$\Theta\subset\C$} iff it is a countable union of irreducible
components of $\Gamma$-invariant co-etale closed subsets defined over
$\Theta$.
\end{defn}
\def\clt{{\Cl_\Theta}}
\begin{defn} For a set $V\subset \U^n$, let $\clt V$ be the
intersection of all closed $\Theta$-definable sets containing $V$:
$$\clt(V)=\bigcap_{V\subset W, W/\Theta{\text{ is }}
\Theta{\text{-definable closed}}} W$$

A point $v\in V$ is called \emph{$\Theta$-generic} iff $V=\clt(v)$,
i.e.~there does not exist a closed $\Theta$-definable proper subset
of $V$ containing $v$.
\end{defn}


\begin{lem}\label{fin.char.cl}\bi\item[(a)] $\clt(V)$ is $\Theta$-definable
\item[(b)]$\clt (V)=\bigcup_{v\in V} \clt(v)=
\bigcup_{S\subset_{\text{fin}} V} \clt(S)$ (union over all finite
subsets)\ei
\end{lem}
\bp $(a):$ By Decomposition Lemma, it is sufficient to consider only
irreducible $V$. However, for irreducible $V$ we may assume that all
sets appearing in the definition of $\Clt(V)$ are again irreducible
and therefore the intersection is finite. It is immediate that a
finite intersection of $\Theta$-definable sets is
$\Theta$-definable.

$(b):$ This follows from the Decomposition Lemma. If $V$ is
irreducible, then $V=\clt(v)$ for $v$ a $\Theta$-generic point of
$V$. If not, by Decomposition Lemma, $V$ decomposes  as a union of
translates of irreducible sets $V_1,\ldots,V_n$. Thus the union
$\bigcup_{v\in V} \clt(v)$ is  the union of the corresponding translates
of the closures $\clt(V_1)$,\ldots,$\clt(V_n)$ of the irreducible components
$V_1,\ldots,V_n$. By  Lemma~\ref{lem:topology}, $\clt(V_i)$ being
closed implies any union of translates of $\clt(V_i)$ is closed; and
thus $\bigcup_{v\in V} \clt(v)$ is a finite union of closed sets,
therefore closed itself. But obviously $V\subset \bigcup_{v\in V}
\clt(v)$ and therefore $\clt(V)\subset\bigcup_{v\in V} \clt(v)$. On
the other hand, for any $v\in V$ $\clt(v)\subset\clt(V)$, and thus
$\clt(V)\supset \bigcup_{v\in V} \clt(v)$. This implies the
lemma.\ep

\begin{lem}\label{def:barq} If a set $W'\subset\uU$ is defined over
$\barq\subset\C$ then $W'\subset\uU$ is $\La$-defined with
parameters from $\pinv(\A(\barq))$. \end{lem} \bp An irreducible
component of the preimage of an algebraic variety $W(\C)\subset\AC$
defined over $\barq$ is an irreducible component of the preimage of the
variety
$$
\bigcup\limits_{\sigma\in\gal(\barq/k)} \sigma W(\C)$$ defined over
$k$. In order for the union to be finite, we use that $W$ is defined
over $\barq$, i.e.~over a finite degree subfield of $\barq$. The
relation $\sim_W$ is in $\La(A)$, and $W'$ can be defined by
$x\sim_W a_1 \aand \ldots \aand x \sim_W a_k$, for some set of
$\barq$-rational points $a_1,\ldots,a_k\in W'(\barq)$. \ep


Recall we assume $\Theta$ to be algebraically closed.
\begin{lem}\label{lem:generic.lift} For every finite index subgroup $H\nrmfin\Gamma$, if $W'$ is irreducible co-etale closed,
then
$w'\in W'$ is $\Theta$-generic iff $w=\phh(w')\in W=\phh(W')$ is
$\Theta$-generic in $W$.\end{lem}

\bp The point $w'\in W'$ is not $\Theta$-generic iff there exists a
$\Theta$-defined irreducible set $w'\in V'\subsetneq W'$;
necessarily $\dim V'<\dim W'$ and $\phh(V')\neq \phh(W')$.\ep

We would rather avoid using this corollary due to its non-geometric
character, but unfortunately we do use it.

\begin{lem}\label{lem:1.5.6} A connected component of a non-empty $\Theta$-generic fibre of a
co-etale closed irreducible set defined over $\Theta$ contains a
$\Theta$-generic point. That is, if $W'\subset \U\times \U$ is co-etale 
irreducible and $\pr:W'\ra \U$ is the projection, and $g'\in \Cl\pr
W'$ is a $\Theta$-generic point of the co-etale closed set $V'=\Cl\pr
W'$, then the {\em $\Theta$-generic fibre} $W'_{g'}=\pr\inv(g')$
contains a $\Theta$-generic point of $W'$.
\end{lem} \bp Basic properties of generic points of algebraic varieties imply 
this property for algebraic varieties.
Let $\wgdc$
be a connected component of a fibre of $W'$ over a $\Theta$-generic
point $g'$ of $\Cl\pr W'$. Then $p(\wgdc)$ is a connected component
of the fibre $\wg$, where $W=\phh(W'),g=p(g')$ is such that $W'$ is
a connected component of $\phh\inv(W)$; this may be seen with the
help of the  path-lifting property, for example. Genericity of
$g'\in\Cl\pr W'$ implies that the point $g\in \Cl\pr W$ is
$\Theta$-generic, and, as a connected component of the fibre $\wg$
of an algebraic variety, $p(\wgdc)$ contains a $\Theta$-generic
point, and then its preimage in $\wgdc$ is also $\Theta$-generic.
\ep


\subsection{(WP) Weak properness : Stein factorisation and fundamental groups }


Above establishes that $\U$ satisfies all but those axioms of an analytic 
Zariski structure that describe the image of a projection --- (SP),(WP) 
and (FC). To prove these these axioms, we use that in algebraic geometry, {\em all} 
morphisms are {\em topologically}
very {\em simple}: {\em each} morphism of complex smooth connected algebraic varieties is,
excepting a closed subset of smaller dimension,  
a topological {\em fibre bundle} with {\em connected} fibres, followed by 
a finite {\em topological covering} (i.e., a fibre bundle with finite fibres).
This is known as {\em Stein factorisation}. Via the long exact  sequence of  
a fibration, this allows us to describe the behaviour of the fundamental group with respect to 
algebraic morphisms. We use this to prove (FC).


Let us give an idea behind the calculations. We need to exclude the counterexample of a finite non-closed spiral 
in $\Cm\times\Cm$ projecting onto a circle in $\Cm$. In the cover, the spiral $S$ unwinds to a 
curve $S'$ of finite length while the circle $S^1$ unwinds to an infinite line $L$. 
As countably many deck translates of $\pr S'$ cover the whole of the line $L$, 
their dimension must be the same in an analytic Zariski structure. Observe that for the counterexample it is 
essential that the projection $\pr \pi(S) \longrightarrow  \pi(S')$ is not surjective, a possibility excluded by
Proposition~\ref{lem:MainGroupProperty}.  

Let us remark that although the circle is not definable for obvious reasons, the variety $\Cm$ is definable  and homotopic 
to the circle, and so considerations above imply that  we need to show 
there is no irreducible co-etale closed subset of $\C^n$ with finite
deck transformation group projecting surjectively onto $\Cm$.


\subsubsection{Prerequisites: topological structure of algebraic morphisms}

\paragraph*{ Exact sequence of fundamental groups of a fibration.}

\begin{fact} \label{fact:short-exact-sequence} For a fibration $f:A \longrightarrow  B$ of 
(nice, e.g. Hausdorff, linearly connected, locally linearly connected) topological spaces,
a pair of points $a\in A,b=f(a)\in B$,
we have an exact sequence of homotopy groups 
$$
\longrightarrow \pi_2(B,b)\longrightarrow \pii(f\inv(b),a)\longrightarrow \pii(A,a)\longrightarrow \pi(B,b)\longrightarrow \pi_0(f\inv(b),a)\longrightarrow \pi_0(A,a)\longrightarrow \pi_0(B,b)\longrightarrow 0
$$
\end{fact}

\begin{rem} In fact, fibrations are thought of as  analogues of exact sequences of Abelian groups
in 'the non-Abelian context' of topological spaces. 
\end{rem}
 
\paragraph*{ { Normal closed analytic sets\label{var:normal}}.}
\begin{defn}(\cite[\S7.2,Def.7.4]{demailly}) A closed analytic subset $X$ is
\emph{normal} at a point $x\in X$  if 
the ring $O_{X,x}$ of germs of holomorphic functions over neighbourhoods of $x$
 is integrally closed in its field
of fractions. A closed analytic subset is normal iff it is normal at every point.

A \emph{normalisation morphism $\nnn$ of variety $Y$} is a morphism
$\nnn:X\ra Y$ from a normal variety $X$ such that any dominant,
i.e. surjective on a Zariski open subset, morphism $f:Z\ra Y$ lifts up to a
unique morphism $\tilde f:Z\ra X$ such that $f=\tilde f\circ \nnn$.
%
\end{defn}

Any smooth closed analytic set is normal (\cite[\S7.4]{demailly}).

We only use the following two properties of a normal variety:

\begin{fact} A normalisation morphism exists for any
variety, and is functorial. Namely, for every variety $(Y,y),y\in Y$ with a base-point 
we may choose a normalisation morphism $\nnn: (\nnn(Y),\nnn(y)) \longrightarrow  (Y,y)$ such that   
for every pair of morphisms $f:(X,x) \longrightarrow  (Y,y), g:(Y,y) \longrightarrow  (Z,z)$ it holds that
$\nnn(fg)=\nnn(f)\nnn(g)$. 
\end{fact}
\bp Lemma~\S7.11 of \cite{demailly} and 
Oka's normalisation principle of [loc.cit.,\S7.12].\ep

%
%


\begin{fact}\label{fact:nrmconn} Let $X$ be a closed analytic subset of a Stein manifold, or
let $X$ be an algebraic variety. If $X$ is connected and normal,
then $X$ is irreducible.
\end{fact}
\bp  Implied by \cite[\S7.4]{demailly}. \ep

\paragraph*{ Fundamental groups of open subsets of normal varieties.}

\begin{fact}\label{fact:piiopen}
Let $Y$ be a connected \emph{normal} complex space and 
$Y^0\subset Y$ be open. Then $\pii(Y^0(\C),y_0) \longrightarrow  \pii(Y(C),y_0)$ is
surjective, for every $y_0\in Y^0(\C)$.
\end{fact}
\bp Kollar, Prop.2.10.1\ep

\paragraph*{ Stein factorisation.}

\begin{fact}\label{fact:stein} Any projective morphism $f:Y\ra X$ of algebraic
varieties admits a
factorisation $f=f_0 \circ f_1$ as a product of a finite morphism
$f_0:Y\ra Y'$ and a morphism $f_1$ with connected fibres.
\end{fact}
\bp \cite[Ch. III, Corollary 11.5]{HartAG} \ep

\paragraph*{ A morphism of normal algebraic varieties is topologically a fibration on an Zariski open subset.}

For normal varieties we have a more precise statement:

\begin{fact}\label{stein}\label{fact:Stein.fact} Let $f:X\ra Y$ be a
morphism of irreducible \emph{normal} algebraic complex varieties 
such that $Y\subset f(X)$.

Then there exist an open subset $Y^0\subset Y$ and $X^0=f\inv (Y^0)$,
and a variety $Z^0$ such that $f$ factories as follows:

$$X^0 \ra^{f^0} Z^0 \ra^{f^{et}} Y^0$$

where \bi \item $Z^0\ra Y^0$ is a finite topological covering in complex topology 
	(i.e. an \etale{} morphism)
\item $X^0\ra Z^0$ is a topological fibre bundle (in complex
topology) with connected fibres
\ei

In particular, 

3. $f:X^0\longrightarrow Y^0$ is a fibration, and its fibres are of boundedly many connected components

4. we have a short exact sequence 
$$
\longrightarrow \pi_2(Y^0(\C),y_0)\longrightarrow \pii(f_{|X^0(\C)}\inv(y_0),x_0)\longrightarrow \pii(X^0(\C),x_0)\longrightarrow \pii(Y^0(\C),y_0)\longrightarrow \pi_0(f_{|X^0(\C)}\inv(y_0),x_0)
$$

\end{fact}
\bp Kollar, Proposition 2.8.1. \ep

Note that while $f^0:X^0\ra Y^0$ is interpretable in the theory of
algebraic varieties and in $\la$, as indeed any morphism of
algebraic varieties is, the theory may not say anything about the
induced morphism $(f^0)_*:\univ({X^0})\ra \univ ({Y^0})$ of the
universal covering spaces of $X^0(\C)$ and $Y^0(\C)$.

\paragraph*{ Morphisms of fundamental groups of normal varieties.}

The Fact~\ref{fact:Stein.fact} above leads to a fact about fundamental groups
 specific to algebraic geometry.

\begin{fact}\label{fact:exact}
Let $f: X\ra Y$ be a morphism of normal algebraic connected complex
varieties; 
assume that $f(X)$ is open is $Y$.

Then there is an open subset $Y^0\subset Y$ defined over the same field as $Y$,
such that for every point $g\in Y^0(\C) \subset Y(\C)$, every point
$g' \in X_g=f\inv (g)$  a generic fibre of $f$ over generic  point $g\in Y(\C)$,
it holds that  
the sequence
$$
f_*:\pii(X_g(\C),g')\ra \pii(X(\C),g')\ra \pii(Y(\C),g) \ra 0 $$
is exact up to finite index.
\end{fact}

\bp Follows from Facts~\ref{fact:piiopen} and~\ref{fact:Stein.fact}
and~\ref{fact:short-exact-sequence}(the exact sequence of the fundamental groups of a fibration).
That is,  Kollar, Proposition 2.8.1 and Kollar, Proposition 2.10.1. \ep

\subsubsection{\label{sec:nonnormal} Extending to non-normal subvarieties }

The above provides an explicit description of morphisms topologically,
between normal algebraic varieties.

However, we need  to deal with an
\emph{arbitrary} subvarieties, not necessarily normal. We do so by
considering the image of the fundamental groups in the big ambient
variety that is normal.

\paragraph*{ Fundamental subgroups of non-normal subvarieties.}

\begin{fact} \label{fact:sobaki} Assume $\A$ is LERF.

Let $p:\U\ra \A(\C) $ be the universal covering space, let $\ita:W\ra
\A\times \A$ be a closed subvariety, and let $Z=\Cl\pr W$. Assume that
 $\p\inv(W(\C))$ and $\p\inv(Z(\C))$ are unfurled.
 

Then there is an open subset $Z^0\subset Z$ 
 defined over the same field as $Z$,
such that for every point $g\in Z^0(\C) \subset Z(\C)$, every point
$(g,g') \in W_g=f\inv (g)$  a generic fibre of $f$ over generic  point $g\in Z(\C)$,
it holds that
the sequence of 
 subgroups of $\pii(\A(\C) ^2,(p(g'),p(g))$
$$
\itastar\pii(W_g(\C),(g,g'))\ra \itastar\pii(W(\C),(g,g'))\ra
\itastar\pii(Z(\C),g)\ra 0
$$
which is exact up to finite index, and the homomorphisms are those
of subgroups of $\pii(\A(\C) ^2,(p(g'),p(g))$.

%
%
\end{fact}
\bp
We prove this by passing to the normalisation of varieties $W$ and
$Z=\Cl\pr W$. The assumption about the irreducibility of connected
components implies that the composite maps of fundamental groups
$\pii(\hatW)\ra \pi(W)\ra \ipii(W)$ and $\pii(\hatZ)\ra \pii(Z)\ra
\ipii(Z)$ are surjective.

To show this, first note that the universal covering spaces $\tilde
\hatW(\C)$ and $\tilde \hatZ(\C)$ are irreducible as analytic
spaces; indeed, normality is a local property, and so they are
normal as analytic spaces; they are obviously connected, and for
normal analytic spaces connectivity implies irreducibility.

By properties of covering maps, a morphism between analytic spaces
lifts up to a morphism between their universal covering spaces (as
analytic spaces); thus the normalisation map $\nnn_W:\hatW \ra W$
lifts up to a morphism $\et\nnn_W:\et\hatW\ra \U$. The normalisation
morphism $\nnn_W$ is finite and closed by Hartshorne
\cite[Ch.II,\S3,Ex.3.5,3.8]{HartAG}; therefore $\et\nnn_W$ is also,
and the image of an irreducible set is irreducible. Therefore
$\et\nnn_W(\et\hatW)$ is an irreducible subset of a connected
component of $\pinv(W(\C))$. Moreover, if we choose different
liftings $\et\nnn_W$, we may cover $\pinv (W(\C))$ by a countable
number of such sets. Now, we use the assumption that a connected
component of $\pinv(W(\C))$ is irreducible to conclude that the
image $\et\nnn_W(\et\hatW)$ coincides with a connected component of
$\pinv (W(\C))$. This implies that the map of fundamental groups is
surjective; this may be easily seen if one thinks of a fundamental
group as the group of deck transformations.

\def\wg{{W_g}}
Let $\nnn_W:\hatn W\ra W$, $\nnn_\wg:\hatn \wg\ra \wg$
and $\nnn_Z:\hatn Z\ra Z$ be the
normalisation of varieties $W$,$\wg$ and $Z$.

By the universality property 
of normalisation  in \S\ref{var:normal} we
may lift the normalisation morphism $\nnn_\wg:\hatn\wg\ra \wg$ to
construct a commutative diagram:
\begin{center}
$%
\begin{array}{cccccc}
   & \hatn\wg & \ra & \hatW & \ra & \hatZ \\
   & \downarrow &  & \downarrow &  & \downarrow \\
   & \wg & \ra & W & \ra & Z \\
\end{array}%
$
\end{center}

By functoriality of $\pi_1$, this diagram and embedding $\ita:W\ra \A\times \A$
gives us

\begin{center}
$%
\begin{array}{cccccc}
   & \pii(\hatn\wg) & \ra & \pii(\hatW) & \ra & \pii(\hatZ) \\
   & \downarrow &  & \downarrow &  & \downarrow \\
   & \pii(\wg) & \ra & \pii(W) & \ra & \pii(Z) \\
   & \downarrow &  & \downarrow &  & \downarrow \\
   & \ipii(\wg) & \ra & \ipii(W) & \ra & \ipii(Z) \\
\end{array}%
$
\end{center}

Now, $g'$ is $\Theta$-generic in $\hatn W'_{g'}$; We are almost
finished now. By Fact~\ref{fact:exact} the upper row of the diagram
is exact up to finite index, and $\pii(\hatW) \ra \pii(\hatZ)$ are
surjective, up to finite index; by assumptions on $W$ and $Z$, the
composite morphisms $\pii(\hatZ)\ra \ipii(Z)$ and $\pii(\hatW)\ra
\ipii(W)$ are surjective. Diagram chasing now proves that the bottom
row is also exact up to finite index, and the map $\ipii(\hatW) \ra
\ipii(\hatZ)$ is surjective up to finite index.
\ep

\subsubsection{Deck transformation groups of co-etale irreducible sets  }

Recall notation $\pi(V')=\{\gamma\in\Gamma^n: \gamma V'\subset V'\}$ for $V'\subset \U^n$, 
and that
if $V'$ is a connected component of $\pinv(V(\C))$, then
the deck transformation group $\pi(V')$ is canonically identified 
with the fundamental group  $\pi_1(V(\C),x_0)$, $x_0\in p(V')$:
to an element $\gamma\in \pi(V')$ there corresponds path $p(\gamma_{x'_0,\gamma x'_0})$
where $x'_0\in V'$ is arbitrary such that $p(x'_0)=x_0$.

\paragraph*{ Deck transformation group of a co-etale irreducible set is cocompact .}

\begin{cor} In a co-etale irreducible set $W$, the deck transformation group $\pi(W)$ 
acts cocompactly, i.e. transitive up-to-compact. 

That is, for every co-etale irreducible closed set $W$ there is a compact subset $W_O\subset W$
such that every point $w\in W$ there are $ \gamma w_O$, $\gamma \in \pi(W), w_O\in W_O$
and $w=\gamma w_O$.
\end{cor}

\bp By Decomposition Lemma, $W$ is a connected component of $H W=\phh\inv\phh(W)$, for
some finite index subgroup $H<\Gamma$. As $\phh$ is a local isomorphism, 
$HW$ being closed analytic implies $\phh(W)=\phh(HW)\subset \U^n/H^n$ 
is closed analytic and therefore Zariski closed by Chow Lemma.
This implies that $W$ is a topological covering of a closed set compact in complex 
topology, and $\pi(W)\cap H$ is its deck transformation group. This implies the corollary.
\ep


\paragraph*{ Deck transformation  group of the projection of an irreducible co-etale closed set.}

\begin{prp}[Action of $\pi(\uU)$ on $\uU$]\label{lem:MainGroupProperty}
\label{prp:MainGroupProperty}\label{lem:mainpi}

Let $W'$ and $V'=\Cl\pr W'$ be  co-etale irreducible closed sets.
Then there is a finite index subgroup $H\nrmfin \Gamma$ such that

\bi\item $\pi(W')\cap H= \{\gamma\in H: \gamma W'\subset
W'\}=\{\gamma\in H: \gamma W'\cap W'\neq \emptyset\}
=\{\gamma\in H: \gamma x'_0\in W'\}$, for any point $x'_0\in W'$

\item $\pr [\pi(W')\cap H] = \pi(V')\cap H$.

\item for an open subset $V^{0'}\subset V'$ it holds that for arbitrary
connected component $\wcgd$ of fibre $\wgd$ over $g'\in V^0$ there
is a sequence exact up to finite index \becd \pi(\wcgd) @>>> \pi(W')
@>\pr_*>>\pi(V')@>>>0,\eecd i.e.~there exists a finite index
subgroup $H\nrmfin \Gamma$ independent of $g$ and $\wcgd$ such that the
sequence is exact: \becd \pi(\wcgd)\cap H @>>> \pi(W')\cap [H\times
H] @>\pr_*>>\pi(V')\cap H@>>>0,\eecd

Moreover, if $W'$ and $V'$ are defined over an algebraically closed
field $\Theta$, so is $V-V^0$. In particular, the above sequence is exact
for $g$ a $\Theta$-generic point of $V'=\Cl\pr W'$.
\ei\end{prp}

\bp[Proof of Proposition] To  prove $(1)$, apply Decomposition Lemma
to  the co-etale closed set $\Gamma W'$; by Decomposition Lemma, take
$H\nrmfin\Gamma$ to be such that the set $\Gamma W'$ decomposes as a union
of a finite number of $H$-invariant sets whose connected components
are irreducible, and therefore they are translates of $W'$. This
implies $(1)$. The item $(2)$ is implied by $(3)$.

Let us now prove item $(3)$.
Let $H$ be such that $W'$ and $V'$ are
connected components of $\phh\inv W(\C)$, $\phh\inv (V(\C))$, respectively,
where $W(\C)=\phh(W'),V(\C)=\phh(V')$. Consider projection morphism
$\pr:\A\times \A\ra A$; it induces a morphism $\pr:W(\C)\ra V(\C)$. By
Lemma~\ref{fact:sobaki} it gives rise to a sequence exact up to
finite index:
$$
\itastar\pii(\wcg(\C),w)\ra \itastar\pii(W(\C),w)\ra
\itastar\pii(V(\C),\pr w)\ra 0
$$
where $\wcg$ is a connected component of a fibre of $W$ over $g\in
V$, and $g$ varies in an open subset $V^0$ of $V$, and $w$ varies in
$\wcg$. The index depends only on the Stein factorisation of the
projection, and is therefore independent of $g$ and fibre $\wgcd$.

Recall 
that there is a canonical identification of $\pi(W')$ and
$\istar\pi_1(W(\C),w)$, and of $\istar\pi_1(\wcg(\C),w)$ and
$\pi(\wgcd')$, etc. As a canonical identification respects morphisms,
Proposition is implied.
\ep 

\begin{cor}\label{cor:stein}\label{cor:stein.nrm}\label{cor.stein.nrm} Let $W'$ be a co-etale  irreducible
closed set, and let $V'=\Cl \pr W'$. Then $\pi(\pr W')$ is a finite
index subgroup of $ \pi (V')$. 
\end{cor}
\bp By item (3) of Lemma~\prettyref{lem:mainpi}.
\ep

\subsubsection{\label{sec:chevalley}Corollary: (WP) Weak Properness, i.e. Chevalley Lemma}

\begin{cor}[Chevalley Lemma]\label{lem:chevalley} For the co-etale topology , it holds:
\bi \item[(SP)] Projections of closed irreducible sets are irreducible closed.
\item[$(SP)_{alg}$] Projections of closed sets invariant under a finite index subgroup of the fundamental group, are closed

\item[$(SP)_{gen}$] Projection of an irreducible constructible set contains all generic points of the projection.

\item[(WP)] The projection of an irreducible set open in its closure contains an open subset of 
the closure of the projection. 
closure.

 \ei\end{cor}

\bp  It is easy to check that the projection of an $H$-invariant closed set 
is closed; indeed, say for $H=\Gamma$,  note $\pr p(\Gamma W')=p \pr
(W')$, and thus $\pr \Gamma W'=\pinv p( \pr W')=\pinv p (V) $, where
$V=\pr p(W')$. As $\AC$ is projective, $V$ is a closed algebraic
subset of $\AC$, and thus $\pinv p(V)$ is a $\Gamma$-invariant closed
subset of $\uU$. By definition of $\Tt$, it is co-etale closed.
This proves $(SP)_{gen}$.

To prove (SP), let  $W'$ be a co-etale irreducible closed set which is a
connected component of $HW'$. Let $V'$ be the closure of $\pr W'$;
we intend to apply item (3) of Proposition above.

The set $\pr HW'$ is closed, and thus $V'\subset \pr HW'$. The set
$V'$ is closed, and thus it is contained in a connected component
$V'_1$ of $\pr HW'$.

Take $v'\in V'\subset V'_1$, and find  $w'\in W'$ such that $\pr
(hw') =v'$; this is possible due to $V'\subset \pr HW'$. Also $\pr
W'\subset V'$, and thus $\pr (w')\in V', \pr(h) \pr(w')=v'\in V'.$
Then $v'\in \pr(h)V'_1\cap V'_1$. We may further take $H$
sufficiently small so that $$\pi(V'_1)\cap H=\{\tau\in\Gamma:\tau
(V'_1)\cap V'_1\neq \emptyset\}=\{\tau\in\Gamma:\tau V'_1=V'_1\}.$$
Then $\pr(h)\in\pi(V'_1)$, and
Proposition~\ref{lem:MainGroupProperty}(2) implies there exists an
element $h_1\in \Gamma (W')\cap [H\times H]$ such that $\pr(h)=\pr h_1$.
Then, $h_1 W'=W'$, and thus $\pr (h_1 w')= \pr (h) \pr w'=v'$, as
required.

This argument can be given topologically. We reprove $(SP)_{alg}$ topologically.

  First, we may assume that $W'$ is a connected component of
  $\phh\inv\phh(W')=HW'$, and
  by Chevalley Lemma for algebraic varieties there is a set
   $V^0\subset \pr \phh (W') \subset V$ such that $V^0\subsetneq V$
  is open in $V$. Let $V'$ be the connected component of $\phh\inv(V)$
  containing $\pr W'$.
  Take $V^{0'}=V'\cap \phh\inv (V^0)$; then
  $V^{0'}\subset V'$ is open in $V'$ as an intersection with an open set.

  Take $v'\in V^{0'}$, and take $w'\in W'$,
  $\pr\phh(w')=\phh (v')\in V^0\subset \pr W$; such a point $w'$ in $W'$
  exists by what we call the covering property of connected components.
  Now, $\pr w'\in V'$, and thus
  $\gamma_0\in \pi(V')$ where $\gamma_0$ is defined by
  $v'=\gamma_0\pr w'$. Condition $\pr\phh(w')=\phh (v')\in
  \etAh(K)$ implies $\gamma_0\in H$. Thus the inclusion
  $\pr\pi(W')\cap H=\pi(V')\cap H$ implies there exists
  $\gamma_1\in\pi(W')$, $\pr \gamma_1=\gamma_0$, and thus
  $v'=\gamma_0 \pr w' =\pr (\gamma_1w')$, and the Chevalley
  lemma is proven.

$(SP)_{gen}$ is implied by (SP), as the projection is irreducible and every fibre 
above a generic point of $\pr W$ contains a generic point of $W=\cl S$
that is necessarily contained in $S$. 

(WP) is also implied by (SP). Let $W\subset \U^n$ be irreducible, and 
let $W_i=W_i\cap W\subset W$ 
be closed irreducible subsets of 
$W$ such that $\bigcup_i W_i$ is closed. We need to prove that 
$\pr (W\setminus \bigcup_i W_i)\subset \U^m$ is open in its closure. 
It is easy to notice that
that we may assume 
that $\bigcup_i W_i$ is $\Gamma_g=\ker(\pi(W) \longrightarrow  \pi(\pr(W))$-invariant:  
 use $\gamma$-invariance of every fibre $W_x=\pr\inv(x)$ to check that 
the projection 
$$\pr (W\setminus \bigcup_i W_i) = \pr (W\setminus  \bigcap_{ \gamma\in\Gamma_g} (\bigcup_i  \gamma W_i))$$
does not change: if $W_x=(\bigcup_i W_i)\cap W_x$ then 
$$W_x=\bigcap_{\gamma\in\Gamma_g}\gamma((\bigcup_i W_i)\cap W_x)=\bigcap_{\gamma\in\Gamma_g}(\bigcup_i \gamma W_i)\cap \bigcap_{\gamma\in\Gamma_g}\gamma W_x=\bigcap_{\gamma\in\Gamma_g}(\bigcup_i \gamma W_i)\cap W_x.$$
The infinite intersection is closed in co-etale topology, and therefore every compact 
subset of $\U^n$ intersects only finitely many of closed subsets of 
$\bigcap_{\gamma\in\Gamma_g} (\bigcup_i \gamma W_i)$. 

Take an open ball  $B\subset \U^m$ such that its closure is compact, 
and  take a finite index subgroup
$H<\Gamma$ such that for every $W_i$ intersecting $\pr\inv (B)$ and every $\gamma\in H^n$,
either $\gamma W_i=W_i$ or $\gamma W_i\cap W_i=\emptyset$; we may do so by 
Proposition~\ref{lem:mainpi} taking into account that we only need to consider finitely many 
$W_i$'s by $\Gamma_g$-invariance. Finally consider 
the quotient 
$W_H=(W\setminus \bigcup_{\gamma\in\Gamma_g,W_i\cap \pr\inv(B)\neq\emptyset}\gamma W_i)/H^{n-m}$.  
It is a subset of an algebraic variety, and, in Zariski topology, is open in its closure. 
Therefore by [HartAG,Ch.I,Ex.3.10,Ch.II,Ex.3.22]  $\pr (W_H)$ is constructible 
and also
$\pr (W\setminus \cup_i W_i) \cap B = H^m\pr(W_H)\cap B$. Thus,  
for every open ball $B$, inside of $B$ the set 
$V_b=\pr W \setminus \pr(W \setminus \cup_i W_i) $ 
coincides with a union of co-etale constructible sets, and, consequently,  
the closure of $V_b$ in complex topology coincides with a finite union of 
co-etale closed sets, locally in every open ball $B$. This implies $V_b$
is closed analytic, and than, in the analytic irreducible decomposition, every 
irreducible component coincides with an irreducible component of a 
co-etale closed set. By definition (I), this implies $V_b$ is closed in co-etale 
topology, as required.
\ep


\begin{rem} It is not 
sufficient to show that locally in complex topology,
 $V_b$ contains an open subset of $\pr W$, as the following counterexample 
shows. Let us explain the picture. Consider a countably infinite family of lines in $\C^2$ passing though
a point, and take the union of countably many intervals lying on these lines.
Then, the union is not contained in the completement of any closed analytic 
set; and it is easy to ensure that, on every compact subset, the union 
is contained in the completement of only finitely many of these lines, i.e. is 
in the completement of a closed analytic set. 
\end{rem}

\subsubsection{Corollary: (FC) Parametrising fibres of particular dimensions }

The proof of (FC)(min) is quite similar to that of (WP).

\begin{cor}[] (FC). For a locally closed irreducible set 
$S\subset \U^n\times \U^m$ and the projection 
$\pr : \U^n\times \U^m \longrightarrow  \U^m$, it holds
\newline
(FC)(min) there exists an open set V such that $V \cap pr S \neq \emptyset$
and for every $v\in V \cap pr S$,
        $dim(pr^{-1} (v) \cap  S) = min \{ dim(pr^{ -1} (a) \cap  S)\}$
\newline
(FC)(>) The set ${a \in  pr S : dim(pr^{ -1} (a) \cap  S) \geq  k}$ is of the form $T \cap  pr S$
for some constructible $T$. \end{cor}
\bp (FC)(min) Let $W=\cl S$ be the closure of $S$ and let $H<\Gamma$ a finite index subgroup 
as provided by Proposition~\ref{lem:MainGroupProperty}. By [HartAG,Ch.I,Ex.3.10,Ch.II,Ex.3.22],
for an open subset $V^0\subset \pr H^{n+m}W$, for every point $v\in V^0$,
it holds that every irreducible component of fibre $(W/H^{n+m})_v$ of the algebraic morphism 
$\pr: W/H^{n+m} \longrightarrow  \U^m/H^m$ is of dimension 
$e=\dim W - \dim \pr W = \dim (W/H^{n+m}) - \dim  \pr (W/H^{n+m}) $. The latter that 
every irreducible component of $\pr\inv (W/H^{n+m})_v$ is of dimension $e$ unless empty,
and $V^0 \cap \pr W$ is as required.   
 The proof of (FC)(>) is similar. 
\ep 


\subsubsection{Uniformity of generic fibres }


Let $\pizero(W')$ denote the set of irreducible components of $W'$.

\begin{cor}[Generic Fibres]\label{cor:gereric.fibres}\label{lem:lef.acf.ana}
In notation of Proposition above, for a 
point $g'\in V'=\Cl\pr W'$ not contained in some proper $\Theta$-definable
closed subset of $V$, the fibre $\wgd$ has finitely many connected
components, and for any connected component $\wgcd$ of $\wgd$, it
holds
$$W'\,\cap\, g'\times H\wgcd \,=\, g'\times\wgcd,$$
$$W'\,\cap\, g'\times H\wgd \,=\, g'\times\wgd.$$

In particular, the formulae above hold for $g'$ a $\Theta$-generic point of $V$.
\end{cor}

\bp Let $H$ be as in Proposition~\ref{lem:MainGroupProperty}. The
fibre $\wgd$ is the intersection of $\wg$ with a coordinate plane,
and therefore is co-etale closed. By Decomposition Lemma, the fibre
$\wgd$ is a union of $H$-translates of a finite number of
irreducible sets $Z'_1,\ldots,Z'_k$. 

To prove the claim, take $h\in H$ such that $Z'_i, hZ'_i\subset
\wgd$. Then $(\id, h)\in H\times H$, and $(\id,h\inv) W'\cap
W'\supset g'\times Z'_i \neq \emptyset$, and by
Proposition~\ref{lem:MainGroupProperty}(1) this implies $(\id,h\inv)
W'= W'$ and $(\id,h\inv)\in\pi(W')$. However, by
Proposition~\ref{lem:MainGroupProperty}(2) $\pi(\wgcd)\cap
H=\ker(\pr_* (\pi(W') \ra \pr(V'))\cap H)$, and thus  $
h\in\pi(\wgcd),h \wgcd=\wgcd$ for any connected component $\wgcd$ of
fibre $\wgd$.

To prove $W'\cap g'\times H\wgcd = g'\times\wgcd$, take $h\in H$
such that $g'\times h\wgd\cap W\neq \emptyset$. Then $(\id, h)\in
H\times H$ and \becd (\id,h) W'\cap W'\supset g' \times h\wgd\cap
\wg \neq \emptyset,\eecd by
Proposition~\ref{lem:MainGroupProperty}(1) this implies $(\id,h) W'=
W'$ i.e.~$(\id,h)\in\pi(W')$. Now
Proposition~\ref{lem:MainGroupProperty}(2), $\pi(\wgcd)\cap
H=\ker(\pr_* (\pi(W') \ra \pr(V'))\cap H)$ gives $h \wgd=\wgd$,
i.e~$h\in\pi(\wgd)$, as required.

In particular, $W'\cap H\wgd = \wgd$ and $W'\cap \wgcd=g'\times \wcgd$\ep

\section{{Core sets: A language for the co-etale topology: $k$-definable sets\label{sec:lang}}}

So far we have analysed the topology on $\U$ (and its Cartesian powers
$\U^n$'s) whose closed sets are rather easy to understand. Now, to
put the considerations above in a framework of model-theory, we want
to define a \emph{language} able to define closed sets in the
co-etale topology. From an algebraic point of view, that corresponds
to defining an automorphism group of $\U$ with respect to the co-etale
topology. The automorphism group is to be an analogue of a Galois
group.

In the terminology of \cite{ZilberBook}, this corresponds to a choice 
of core closed subsets. Our language is smaller than that: core closed 
subsets are definable with parameters (corresponding to core subsets).

Let us draw an analogy to the action of Galois group on $\barq$ as
an algebraic variety defined over $\Q$ endowed with Zariski
topology. The Galois group may not be defined as the group of
bijections continuous in Zariski topology: for example, all
polynomial maps are continuous in Zariski topology; linear and
affine maps $x\ra ax+b$ are such continuous  bijections.

Thus we distinguish certain $\Q$-definable subsets among Zariski closed
subsets of $\barq^3$, and then define  Galois group as the group of
transformation (of $\barq$) \emph{preserving the distinguished
$\Q$-defined subsets (of $\barq^3$)}; in this case the graphs of
addition and multiplication. It is \emph{then} derived, rather
trivially, that this implies that Galois group acts by
transformation continuous in Zariski topology.

Recall the way this is derived: the $\Q$-definable subsets are given
\emph{names}, in this case addition and multiplication, and then
each closed set (subvariety) is given a name by the equations
defining the set of its points; in fact, in algebraic geometry the
word variety means rather the \emph{name}, the set of equations,
rather that the set of points the equations define.

In order to define a useful automorphism group of the co-etale topology,
we follow the same pattern.

Model theory provides us with means to give precise meaning to the
argument above, and to define mathematically what is it exactly that
we want. In these terms, the distinguished subsets form a
\emph{language}, and the Galois group is the group of automorphisms
of \emph{the structure in that language}. Model theory studies that
group via the study of the structure.

\subsection{Definition of a language $\La$ for universal covers in
the co-etale topology\label{III.3}}\def\k{k}

In this $\S$, it becomes essential that $\A$ is defined over an
algebraic field $\k\subset \barq\subset \C$ embedded in $\C$.


We consider $p:\U\ra \A(\C) $ as a structure in the following language.

\begin{defn}\label{def:tau}

We consider the universal covering space $\pUA$ as a one-sorted
structure $\uU$, in the language $\La$ which has the following
symbols: \bi
\item[] the symbols $\simz$ for $Z$ a closed subvariety of $\AC^n$
defined over number field $k$, and, \item[] the symbols $\simH$, for
each normal subgroup $H\nrmfin \pi(\uU)^n$ of finite index \ei

The symbols are interpreted as follows:

\bi\item[] $x'\simz y'$ $\iff$ points $x'\in\uU^n$ and $y'\in\uU^n$
lie in the same (analytic) irreducible component of the
$\Gamma$-invariant closed analytic set $\pinv(Z(\C))\subset\uU^n$.

\item[] $x'\simH y'$ $\iff$ $\exists \tau\in \pi(\uU)^n :\tau
x'=y'$ and $\tau\in H$.\ei
\end{defn}

 Note that we do not assume $Z$ to be
connected.


As justified by Corollary~\ref{lem:simcsim}, we get a bi-interpretable language
by considering the following predicates instead:

\bi \item[] $x' \sim_{Z,\etAh}^c y'$ iff $x'$ and $y'$ lie in the same
 \emph{connected} component of the preimage $\phh\inv (Z_i(\C))$,
 $Z_i\subset \etAh(\C)^n$  an irreducible component of algebraic
 variety $\pHh\inv(Z(\C))\subset \etAh(\C)^n$.\ei

\begin{cor}\label{lem:simcsim} For every closed $\Gamma$-invariant analytic subset $Z'$ of $\U^n$,
there exist closed analytic subsets $Z_1',\ldots,Z'_n$ invariant
under action of a finite index subgroup $H$ of $\Gamma$, such that
$$x\sim_{Z'} y \iff x\simc_{Z'_1} y \OR x\simc_{Z'_2} y \OR \ldots \OR
x\simc_{Z'_n} y.$$ Consequently, for every closed subvariety $Z$ of
$\A$, there exist subvarieties $Z_1,\ldots,Z_n$ of a finite \etale{}
cover $\etAh$ such that
$$x\sim_Z
y \iff x\simc_{Z_1} y \OR x\simc_{Z_2} y \OR \ldots \OR x\simc_{Z_n}
y.$$
\end{cor}
\bp Take $H$ and $Z'_1,\ldots,Z'_n$ as in Decomposition
Lemma~\ref{lem:noeth:ana}.\ep
Note that the language $\La$ is countable. This is an essential
property, from model-theory point of view; in technical,
down-to-earth terms it is useful to make inductive constructions.

Let us use this opportunity to remind that we use symbols $\sim_Z$
rather abusively to mean ``lie in the same irreducible component
of'' either $\Gamma Z$, $\phh\inv(Z)$, etc.

\subsection{\label{sec:la-def}$\La$-definability of $\pi(\uU)$-action etc}

In the next lemma, a closed set means a co-etale closed set.

\begin{lem}\label{sec:la-def:lem} For any normal finite index subgroup $H\nrmfin \Gamma$ it holds \bi
 \item the relation
   $$\Aff_H(x,y,z,t)=\exists \gamma\in H: \gamma x=y \aand \gamma z=t$$
  is $\La(\emptyset)$-definable
 \item An $H$-invariant closed set is $\La$-definable with parameters.
 \item A connected component of a generic fibre of
  an $\La$-definable irreducible closed set is
  uniformly $\La$-definable; the definition is valid 
  over an open subset of the projection, definable
  over the same set of parameters.
  \item Any co-etale closed irreducible set is a connected component of
  a fibre of an $\La$-definable set.
 \item An irreducible closed set is $\La$-definable.
 \ei
\end{lem}

\bp To prove (1), note that
$$\pinv(\Delta(\C))=\bigcup\limits_{\gamma\in\Gamma}\{(x,\gamma x): x\in
U\}$$ where $\Delta=\{(x,x):x\in A\}$ is an algebraic closed
subvariety defined over $\k$. The connected components $\{(x,\gamma
x): x\in \U\},\gamma\in\Gamma$ are the equivalence classes of
the relation $\sim_{\Delta}$, and thus are definable with parameters.

Evidently $\Aff_\pi(x,y,s,t)$ iff $(x,y)\sim_{\Delta} (s,t)$ lie
in the same connected component of $\pinv (\Delta(\C))\subset
U\times U$.

To prove (2), we consider two cases. {\bf $\barq$-case:} An
irreducible closed subvariety $\factor Z \barq\subset \A$ defined
over $\barq$ is an irreducible component of subvariety
$$Z_{\k}=\bigcup\limits_{\sigma: k_Z\hra\C} \sigma(Z)$$ of $\A$,
where $k_Z$ is the field of definition of $Z$ of finite degree. The
formula implies $Z$ is $\La$-definable with parameters with the help
of symbol $\sim_{Z_\k,A}$; the parameters may be taken to lie in $\A(\barq)$
but not necessarily in $\A(k_Z)$. A slightly more complicated
argument could give a construction defining $Z$ as a connected
component.

For an analytic co-etale closed irreducible set $Z'\subset\uU$, it
holds that $Z'$ is an irreducible component of $\Gamma Z'$, i.e.~it is
an irreducible component of $\pinv(Z)=\pinv p(Z')$. Thus the above argument gives
that every co-etale irreducible subset of $\uU$ defined over $\barq$
is $\La$-definable with parameters.

 {\bf $\overline{\Q(t_1,\ldots,t_n)}$-case:}
    Thus we have to deal with the case when $p(Z)$ is not
    $\barq$-definable. Our strategy is to show that any such set is
    a connected component of a $\barq$-generic fibre of a
    $\barq$-definable set, and then show that such connected
    components are uniformly definable. Uniformity will be important
    for us later in axiomatising $\uU$.

    Let us  see first that each co-etale closed irreducible set is a
    connected component of a fibre of a co-etale closed irreducible set
    defined over $\barq$.

    Take a co-etale irreducible set $Z'$ and take $H\nrmfin\Gamma$ such that $Z'$ is a connected component
    of $HZ'=\phh\inv(Z)$, for an irreducible algebraic closed set $Z=\phh(Z')$.
    By the theory of algebraically closed field, we know that $Z$ can be defined as
    a Boolean combination, necessarily a positive one, of $\barq$-definable closed subsets
    and their fibres; by passing to a smaller subset if necessary,
    we see that the irreducibility of $Z$ implies that
    algebraic subset $Z\subset \\A(\C) $ is a connected component of
    a $\barq$-generic fibre of a  $\barq$-definable closed subset
    $W\subset \\A(\C) ^n$.  Then $HZ'$ is the corresponding fibre of $\phh\inv(W)$.
    The closed set $Z'$ is a union of
    the corresponding fibres of the irreducible components of
    $\phh\inv(W)$, and irreducibility of $Z'$ implies that union is
    necessarily trivial. Thus, we have that $Z'$ is a connected
    component of a fibre of an irreducible co-etale closed set defined over $\barq$.
    We may also ensure that $Z'$ is a connected component of a
    $\barq$-generic fibre of $W'$ by intersecting $W'$ with
    the preimage of an irreducible $\barq$-definable set containing
    $\pr Z',$ and repeating the process if necessary.

  Let us now prove that the connected components of the
    $\barq$-generic fibres of an irreducible $\barq$-definable set
    are $\barq$-definable.

    Let $W'\subset \AC^2$, and let $V'=\Cl\pr W'$ be as in
    Proposition~\ref{prp:MainGroupProperty} and
    Corollary~\ref{cor:gereric.fibres}. The morphism $\pr: W\ra V$ admits
    a Stein factorisation (Fact~\ref{fact:stein}) $\pr=f_0\circ f_1$
    as a composition of a finite morphism $f_0:W\ra V_1$ and
    a morphism with connected fibres $f_1:V_1\ra V$. In particular,
    two points $x_1,x_2\in \wg$ lie in the same connected component
    of fibre $\wg$ iff $f_0(x_1)=f_0(x_2)$.
\def\simc{\sim^c}

   Now set
   \belcd x'\simc_{W_g} y' \iff x'\sim_{W} y' 
   \aand \pr x' = \pr y'
   \aand f_0(\phh(x'))=f_0(\phh(y'))\label{fml:stein}\eelcd
   (here subscript $g$ is a part of the notation, and does not
   denote an element of $\U$).

   In notation of Corollary~\ref{cor:gereric.fibres}, we have

\begin{cor}\label{cl:stein} If  $\pr x'=g'\in V'^0$, then
the formula $x'\simc_{ W_g} y'$ holds iff $x'\simH y'$
   and $x'$ and $y'$ lie in the same connected component of fibre $\wgd$
   of $W$.
   If $W,V$ are $\Theta$-definable, so is $V'^0$.
   The parameters needed to define $\simc_{W_g}$ live in $\U/H$.
\end{cor}

\bp This is a reformulation of the formula
$W'\cap g'\times H\wgcd = g'\times\wgcd$. Indeed,
$\pr x' = \pr y'
   \aand f_0(\phh(x'))=f_0(\phh(y'))$ holds
   iff $x',y'\in g'\times H\wgcd $ for $g'=\pr x'=\pr y'$ and
   some $\wgcd$ a connected component of fibre of $W'$ above $g'$.
   The relation of lying in the same connected component of a fibre
   being translation invariant, we may as well assume
   $x',y'\in W'$ if  $x'\sim_{W} y'\in W'$ lie in the same
   connected component of $W'$. Then the formula means
   that $x',y'$ lie in the same connected component of
   fibre $g'\times \wgd$.

  The claim that the formula holds for $g'\in V'^0$  in an open subset is $\Theta$-definable is a part of
  the conclusion of Corollary~\ref{cor:gereric.fibres}.
\ep

The claim above implies (3); (4) and (3) imply (5) and (2). \ep

\begin{cor} Let $\autla(\U)$ be the group of bijections $\phi:\U \ra
U$ preserving relations $\simz\in\La$; then $\autla(\U)$ acts by
transformations continuous in the co-etale topology .\end{cor} \bp
Immediate by Lemma~\ref{sec:la-def:lem}. \ep

The results above justify thinking of $\autla(\U)$ as a \emph{Galois
group of $\U$}.

\begin{rem} Via identifications $\factor \U H \cong \etAh(\C)$, there is a natural inclusion
of a subgroup of $\autla(\U)$ into $\Aut(\C/\Q)$; what can one say
about the common subgroup of $\autla(\U)$ and $\Aut(\C/\Q)$, or
rather a conjugacy class of such subgroups? Is there any relations
between $\autla(\uU)$ and the Grothendieck's fundamental group
$\hat\pi_1(\A_\Q,0)$?
\end{rem}

\section[Model homogeneity]{Model homogeneity: an analogue of $n$-transitivity of
$\Aut_\la(\U)$-action.\label{sec:whom}}

 Now we want to study the action of $\Aut_\La(\U)$ on $\U$, and
analyse orbits of its action on $\U$ and $\U^n,n>1$. In model theory
one would hope that the aforementioned orbits can be analysed in terms
of the language; in presence of a nice topology possessing a properness property 
(WP) or (SP) we may hope to analyse orbits in terms of closed sets.

The situation when this is possible is called homogeneity;
Property~\ref{prop:whom} below states \emph{model homogeneity} of
$\U$. Model homogeneity says, roughly, that two tuples of points lie
in the same orbit (of the action fixing an algebraically closed
subfield) iff there are no obvious obstructions, i.e.~iff they lie
in the same closed sets (defined over an algebraically closed
subfield which we assume fixed).

\begin{defn} We say that $W$  is a \emph{$\Theta$-constructible set} iff
\bi \item the closure $\Cl W$ is defined over $\Theta$
    \item $W$ contains all $\Theta$-generic points of the irreducible
    components of $\Cl W$.
\ei An \emph{irreducible} constructible set is a set whose closure
is irreducible.
\end{defn}

\begin{defn} We say that $w\in W$ is a $\Theta$-generic point of 
an irreducible constructible set iff $w$ does not lie in a proper 
$\Theta$-definable subset of $W$.  

We say that a property holds for a {\em uniform generic point} of $W$ iff
it holds for every point is some open $\Theta$-definable subset of $W$. 
\end{defn}

\begin{lem}\label{lem:geomwhom}  A projection of an irreducible $\Theta$-constructible set is
$\Theta$-constructible. \end{lem}\bp Let $W\subest U\times U$ be an
irreducible set defined over $\Theta$, and let $W_0$ be the set of
all $\Theta$-generic points of $W$; generally speaking, $W_0$ is not
definable. We need to prove that $\pr W$ is also
$\Theta$-constructible. Let $g$ be a $\Theta$-generic point of the
closure of $\pr W$; we know $g\in \pr W$ by (SP) of Lemma~\ref{lem:chevalley}. By
Lemma~\ref{lem:lef.acf.ana} 
we know that the (non-empty) fibre $W_g$ contains a $\Theta$-generic
point of $W$, and thus $g\in \pr W_0$, as required. \ep

The set of realisations of  a complete quantifier-free syntactic
type $p/\Theta$ with parameter set $\Theta$ is
$\Theta$-constructible; and conversely, every $\Theta$-constructible
set can be represented in this form.

Thus, the above lemma is equivalent to $\omega$-homogeneity for such
types.

\begin{defn}\label{def:whom} We say that $\U$ is \emph{homogeneous for irreducible closed
sets over $\Theta$}, or {\em homogeneous for syntactic
quantifier-free complete types over $\Theta$}, or \emph{model
homogeneous} iff either of the following equivalent conditions holds
\bi \item the projection of an irreducible $\Theta$-constructible set is
$\Theta$-constructible; \item for any tuples $a,b\in \U^n$ and $c\in
U^m$ if $\qftp(a/\Theta)=\qftp(b/\Theta)$ then there exists $d\in
U^m$ such that $\qftp(a,c/\Theta)=\qftp(b,d/\Theta)$\ei
\end{defn}

To see that the conditions are equivalent, note that the set of
realisations of a complete quantifier-free type $\qftp(a,c/\Theta)$
is
 $\Theta$-constructible;
its projection contains $a$ and also is $\Theta$-constructible; $a$
is its $\Theta$-generic point; then $\tp(a/\Theta)=\tp(b/\Theta)$
implies $b$ is also $\Theta$-generic, i.e. belongs to the projection.

The above proves the following result.

\begin{prop}\label{prop:whom}  The standard model $p:\U\ra \A(\C) $
in language $\La$ is model homogeneous, i.e.~it is $\omega$-homogeneous for
closed sets over arbitrary algebraically closed subfield
$\Theta\subset \C$.
\end{prop}

\bp Follows directly from Def.~\ref{def:whom} and
Lemma~\ref{lem:geomwhom}.\ep

\begin{cor} The set of realisations  of a quantifier-free type
$\qftp(x/\Theta)$ over $\pinv(A(\Theta))$ consists of
$\Theta$-generic points of some co-etale irreducible closed subset of
$\uU$.\end{cor}\bp Follows from the previous statements.\ep

\section[An $\Lww$-axiomatisation $\Ax(\AC)$]{\label{III.5}An $\Lww$-axiomatisation $\Ax(\AC)$ and stability of the corresponding  $\lww$-class.}

In this $\S$ we introduce an axiomatisation $\Ax(\AC)$ for
$\Lww(\La)$-class which contains the standard
model $p:\U\ra \A(\C) $, and is stable over models and all
models in it are model homogeneous. We then show that the class of
models satisfies $(\aspir{2_{\aleph_0\ra\aleph_1}})$ of Theorem~\ref{Th:1}.

\subsection{Algebraic $\La(G)$-structures\label{sec:finitaryLareducts}}

We know that $\uU/G=\etAg(\C)$ carries the structure of an algebraic
variety over field $\C$. The covering $\etAg(\C)\ra \AC$ carries a
structure in a reduct $\La(G)$ of language $\La$. In fact, similar
interpretation works for an arbitrary algebraically closed field $K$
instead of $K=\C$.

For every finite index subgroup $G\nrmfin \Gamma$, there is a
well-defined covering $\etAg\ra \A$ of finite degree. The space
$\AC$ is projective, and thus $\etAg(\C)$ is also a complex
projective manifold. By Fact~\ref{GAGA}, $\etAg$ has the structure of
an algebraic variety.

Recall that we use the following fact as the defining property of an
\'etale covering: the morphism $B(K)\ra A(K)$ of varieties over an
algebraically closed field $K$ of char 0 is \emph{\'etale} iff there
exists an embedding $i:K'\ra\C$ of the field $K'$ of definition of
$\A$ and $B$ into $\C$ such that the corresponding morphism $i(B)(\C)
\ra i(A)(\C)$ is a covering of topological spaces.

\begin{defn}[Finitary reducts of $\La$]
Let $\pG:\etAg(K)\ra \A(K)$ be a finite \'etale morphism. Let
$\La(G)\subset\La$ be the language consisting of all predicates of
$\La$ of form $\sim_Z$ and symbols $\simH$ for $G\subset H$. Then
$\etAg(K)\ra \A(K)$  carries an $\La(G)$-structure as follows: \bi
  \item $x'\sim_Z y' \iff$ points $x',y'\in \etAg(K)^n$ lie in the same
   irreducible component of algebraic closed subset $\pG\inv(Z(K))$ of
   $\etAg(K)^n$.
  \item $x'\simH y' \iff $ there exist an algebraic morphism
      $\tau:\etAg\ra\etAg$ and a co-etale covering morphism
      $q:\etAg\ra \etAh$ such that  $\tau(x')=y'$ and $\tau\circ q=q$:
      \becd
      \etAg @>\tau>> \etAg \\
      @VV q\text{ \etale{} cover} V        @VV q V\\
      \etAh  @>\id>> \etAh\\
      \eecd
  \ei
\end{defn}

For $G=e$ the trivial group and $K=\C$, the construction above
would degenerate into the interpretation of $\uU\ra \A$ if it were
well-defined.

For $G=\Gamma$, $\etAg=\A$, and thus $\La(\Gamma)$ is just a form of the
language for the algebraic variety $\A$; here the point is that we
have predicates for the relations for irreducible components of
$\k$-definable closed subsets only.

In general, the above is simply a variation of an ACF structure on
$\A$. In particular, all Zariski closed subsets of $(\etAg)^n(K)$
are $\La(G)$-definable.

\subsection[Axiomatisation $\Ax(\AC)$]{Axiomatisation $\Ax(\AC)$ of the universal covering space $\uU$}

We define the axiomatisation $\Ax=\Ax(\AC)$ to be an
$\Lww(\La)$-sentence corresponding to Axiom~\ref{ax:0} and
Axioms~\ref{ax:1.1}-\ref{ax:2.3} below.

\subsubsection{Basic Axioms}

These axiom describe quotations $\U/\simH$ for $H\nrmfin \Gamma$, and
some properties of $\U\ra\U/\simH$.

\begin{ax}\label{ax:0} All first-order statements valid in $\uU$ and expressible in terms of
$\La$-interpretable relations $$ x'\sim_{Z,\etAg} y' := \exists x'' \exists y''
(x''\simZ y'' \aand x''\simG x' \aand y''\simG y'),G\nrmfin \Gamma$$ and
$\simG,G\nrmfin\Gamma$.
\end{ax}

Essentially, these axioms describe $\factor\uU G$ as an algebraic variety.

\subsubsection{Path-lifting Property Axiom, or the covering property Axiom}

\begin{ax}[Path-lifting Property for $W$; Covering Property for $W$]\label{ax:1.1} For
every $\La$-predicate $\sim_W$ and all $G\nrmfin\Gamma$ small enough,
we have an axiom \becd
  x'\sim_{W,\etAg} y' \lra \exists y'' (y''\sim_G y' \aand
                                        x' \sim_W y'')
\eecd\end{ax}

We also have a stronger axiom for \emph{fibres} of $W$; here we use
that the relation ``to lie in the same connected component of a
fibre of a variety'' is algebraic and therefore the corresponding
$G$-invariant relation is $\La$-definable.

\begin{ax}[Lifting Property for fibres]\label{ref:liftfibre}\label{ax:1.2}
For all $G\nrmfin\Gamma$ sufficiently small, we have an axiom \becd
  (x'_0,x_1') \simc_{W_g,\etAg} (y_0',y_1') \lra \exists y_1''
                                   [y_0'\sim_G x_0' \aand
                                    y_1''\sim_G y_1' \aand
                                   (x'_0,x_1') \sim_{W}
                                   (x_0',y_1'')]
\eecd
in a slightly different notation
\becd
  x' \simc_{W_g,\etAg} y' \lra \exists y'' (y''\sim_G y' \aand
                                    \pr x' = \pr y'' \aand
                                     x' \sim_W y'')\eecd\end{ax}

The relation $x' \simc_{W_g,\etAg} y'$ is defined by the
formula (\ref{fml:stein}) (cf.~Claim~\ref{cl:stein}).

\begin{ax}[Fundamental group is residually finite]\label{ax:2.1}
\becd \forall x'\forall y' ( x'=y' \iff \AND\limits_{H\nrmfin \Gamma} x'\simH y')\eecd
\end{ax}

Thus, it says that two elements of $\U$ separated by an element of $H$ for every
$H\nrmfin\Gamma$, have to be equal.

The next property is strengthening  of the previous one; namely, if
an element $b$ is $\simH$-equivalent to an element of a group
generated by $a_1,\ldots,a_n$, then it is actually in the group. In
terms of paths, this has the following interpretation: take loops
$\gamma_1,\ldots,\gamma_n$ and a loop $\lambda$. If for every
$H\nrmfin\Gamma$ it holds that $\lambda$ is $\simH$-equivalent to some
concatenation of paths $\gamma_1,\ldots,\gamma_n$, then it is
actually a concatenation of these paths.

\begin{ax}[``Translations have finite length'', subgroup separability]\label{ax:fltransl}\label{ax:2.2}
For all $N\in \N$ we have an $\lww$-axiom
\becd \forall b \forall a_1 \ldots\forall a_N.  \\
\AND\limits_{H\nrmfin\Gamma}  \bigvee\limits_{n\in\N} \exists h_1
\ldots h_n \left( b\simH h_n \aand h_1=a_1 \aand \AND\limits_{1\leq
i\leq n}\bigvee\limits_{1\leq j < N}(h_i,h_{i+1})\sim_{\Delta}
(a_j,a_{j+1})\right) \\ \ \ \lra  \bigvee\limits_{n\in\N} \exists
h_1 \ldots h_n \left( b = h_n \aand h_1=a_1 \aand \AND\limits_{1\leq
i\leq n}\bigvee\limits_{1\leq j < N}(h_i,h_{i+1})\sim_{\Delta}
(a_j,a_{j+1})\right)
 \eecd
\end{ax}

The next axiom is needed to apply the axioms above. It reflects the fact
that the fundamental groups of varieties are finitely generated, 
a fact we used and prove in the proof of Lemma~\ref{lem:noeth:ana}.
recall that this was proved as a corollary of the fact that
topologically an algebraic variety can be triangulated into finitely many contractible
pieces nicely glued together. 

\begin{ax}[Groups $\pi(W_g)$ are finitely generated]\label{ax:2.3}
For every symbol $\sim_W$ and for each $H\subset\Gamma$ small enough we have an
$\lww$-axiom:
\becd \bigvee_{N\in\N}\exists a_1 \ldots\exists a_N \forall b. \\
\AND\limits_{1\leq i\neq j\leq N} (a_i\simW a_j \aand a_i\simH a_j
\aand \pr a_i=\pr a_j) \aand
 \left( \AND\limits_{i=1}^N
 (b\sim_W a_i \aand \pr b=\pr a_i)
\lra \right. \\
\bigvee\limits_{n\in\N} \exists h_1 \ldots h_n \left.\left( b = h_n
\aand h_1=a_1 \aand \AND\limits_{j=1}^{N}
\bigvee\limits_{j=1}^{N-1}
(h_i,h_{i+1})\sim_{\Delta}
(a_j,a_{j+1}) \aand \pr h_i=\pr h_{i+1}\right)\right) \eecd
\end{ax}
In fact, we may combine the two axioms above into one weaker axiom
which would require {\em subgroup separability with respect to the
subgroups $\pi(W)$}.
%
\subsubsection{Standard model $\uU$ is a model of $\Ax$}

The universal covering space $p:\uU\ra \AC$ satisfies the Axiom~\ref{ax:0} by
definition.

To prove $\uU$ satisfies Axiom~\ref{ax:1.1}, note that for $G\nrmfin
\Gamma$ small enough, the relations $x' \sim_{W,G} y'$ means that
$\pgg(x')$ and $\pgg(y')$ lie in the same irreducible component
$W_i$ of the preimage of $W\subset\AC^n$ in $\etAg(\C)^n$. Take a
path $\gamma$ connecting $\gamma(0)=\pgg(x')$ and
$\gamma(1)=\pgg(y')$ lying in $W_i$; by the lifting property it
lifts to a path $\gamma', \gamma'(0)=x'$ such that
$\pgg(\gamma'(t))=\gamma(t)$, $0\leq t\leq 1$. Then,
$\pgg(\gamma'(1))=\pgg(y')$, and thus $\gamma'(1)\simG y'$. On the
other hand, $\gamma'(1)$ and $x'$ lie in the same connected
component of the preimage of the irreducible component $W_i$ in
$\uU$. Now note that by Decomposition Lemma~\ref{lem:decomp.ana} for
$G$ small enough such a connected component has to be irreducible,
and thus Axiom~\ref{ax:1.1} holds.

The Axiom~\ref{ax:1.2} has a similar geometric meaning as Axiom~\ref{ax:1.1}; the assumption
is that $\pgg(x')$ and $\pgg(y')$ lie in the same connected component of a fibre $W_g$;
it is enough to take $\gamma$ to lie in fibre $W_g$ to arrive to the conclusion of Axiom~\ref{ax:1.1}.

Axiom~\ref{ax:2.1} follows from the condition 2 of the definition of a 
LERF variety.

Axioms~\ref{ax:2.2} is condition 2 of the definition of a LERF variety.

The geometric meaning of $(h_i,h_{i+1})\sim_\Delta (a_i,a_{i+1})$ is
as follows. The pair of points $a_i,a_\ii$ determines a path
$\gamma$ in $\AC$, $\gamma(0)=\gamma(1)=p(a_i)=p(a_\ii)$. For points
$h_i,h_\ii$ such that $p(h_i)=p(h_\ii)$, they can be joined by a
lifting of $\gamma$ iff $(h_i,h_{i+1})\sim_\Delta (a_i,a_{i+1})$.
\ldots Thus the assumption in the axiom says that if any two points
of fibre above $p(b)=p(a_1)$ can be joined by a concatenation of
liftings of finitely many paths $\gamma_i$'s in $\AC$, \emph{up to a
translate by an element of $H$}, then they can in fact be just
joined by such a sequence. In a way, this can be thought of as
disallowing paths of infinite length.

On the other hand, the condition $(h_i,h_{i+1})\sim_\Delta (a_i,a_{i+1})$
can be interpreted as $h_\ii=\tau_i h_i$ where $\tau_i$ is the deck transformation taking
$a_i$ into $a_\ii$, $\tau_ia_i=a_\ii$. Then, the assumption says that if $b\in\pi(\uU)$ belongs
to the group generated by $\tau_i$'s, up to $\simH$, then $b$ does belong
to the subgroup generated by $\tau_i$'s.

The last remaining Axiom~\ref{ax:2.3} means that the fundamental groups $\pi(W_g)$
is finitely generated, and we already used this Fact in the proof of Lemma~\ref{lem:noeth:ana}.

\subsection{Analysis of models of $\Ax$\label{ss:Xi.fin}}

\subsubsection{Models  $\factor \U{\,\simH}$ as algebraic varieties}

Let $\U\models \Ax$ be an $\La$-structure modelling axiomatisation
$\Ax(\AC)$, and let $\uU$ be the standard model, i.e.~the universal
covering space of $\\A(\C) $ considered as an $\La$-structure.

We know that $\factor{\uU}{\,\sim_H} \cong \etAh(\C)$  for some
algebraic varieties $\etAh(\C)$ defined over $\C$. The relations
$\sim_H,\sim_{Z,H}$ are essentially relations on
$\factor{\uU}{\,\sim_H}$, and thus Axiom~\ref{ax:0} says that the
first-order theories of $\factor{\uU}{\,\sim_H}$ and that of standard model
$\factor{\uU}{\,\sim_H}$ in the language
$\La(H)=\{\sim_H,\sim_{Z,H}:\text{ Z varies}\}$ coincide. We know by
properties of analytic covering maps that an irreducible co-etale
closed subset of $\uU$ covers an irreducible Zariski closed subset
of $\etAh(\C)$, and thus the relation $\sim_{Z,H}$ on
$\factor{\uU}{\,\sim_H}$ interprets as saying that $x,y\in \etAh(K)$ lie in
the same (Zariski) irreducible component of the preimage of $Z(K)$
in $\etAh(K)$. In particular, every component is definable by $g \simZ y$ where $g$ 
is taken to be its generic point. Since every $\barq$-definable closed subvariety is 
an irreducible component of a $\Q$-definable subvariety, this implies that 
every $\barq$-definable closed subvariety of $\etAh(\C)$ is $\La(H)$-definable. 
Thus, full theory of an algebraically closed field is
reconstructible in $\La(H)$ on $\U/{\simH}$; and thus, there is an
algebraically closed field $K=\bar K, \textrm{char} K=0$ such that
$\U/\simH\cong \etAh(K)$.

Fix these isomorphisms $\U/\simH\cong \etAh(K)$, and let $\phh:\U \ra
\etAh(K)$ be the projection morphism. Then the above considerations
say \bi

\item[] $ x'\sim_{W,H} y' \iff \phh(x')\sim_{W,H} \phh(y')
\iff$ $x'$ and $y'$ lie the same (Zariski) irreducible component of
the preimage of $Z(K)$ in $\etAh(K)$.
\item[] $x \simG y' \iff$ there exist an algebraic morphism
      $\tau:\etAg\ra\etAh$ and a co-etale covering morphism
      $q:\etAh\ra \etAg$ such that  $\tau(x')=y'$ and $\tau\circ q=q$:
      \becd
      \etAh @>\tau>> \etAh \\
      @VV q\text{ \etale{} cover} V        @VV q V\\
      \etAg  @>\id>> \etAg\\
      \eecd
  \ei

An important corollary of above considerations is that any set of
form $\phh\inv(Z(K)), Z(K)\subset\etAh(K)$ is $\La$-definable.

\begin{notation} Let us introduce new relations on $\U$; eventually
we will prove that they are first-order definable. We introduce the
relations below for every closed subvariety of $\A(K)$, not
necessarily defined over $\k$ (those would be in $\La$) \bi
\item[]
  $x'\simW y' \iff \phh(x')\sim_{W,H} \phh(y')$ for all
  $H\nrmfin\Gamma$.\ei
\end{notation}

An \emph{irreducible component} of relation $\simW$ is a maximal set
of points in $\U$ pairwise $\simW$-related. A subset of $\U$ is
\emph{basic closed} iff it is a union of irreducible components of
relations $\sim_{W_1},\ldots,\sim_{W_n}$, for some $W_1,\ldots,W_n$.
An \emph{irreducible closed set} is an irreducible component of a
relation $\simW$ for some closed subvariety $W$. Let us call a
subset of $\U$ \emph{co-etale closed} iff it is the intersection of 
basic closed sets. 
 This defines an analogue of the
co-etale topology  on $\uU$.

\subsubsection{Group action of fibres of $p:\U\ra \A(K)$ on $\U$}

For a point $x_0\in \U$, let $\pi(\U,x'_0)=\{y:y\sim_\Gamma x'_0\}=\pinv
p(x'_0)$ be the fibre of $p:\U\ra \A(K)$. For every point $z'\in \U$
and every point $y' \sim_\Gamma x'_0$, there exists a point $z''\in \U$
such that $\pgg(z',z')\sim_\Delta\pgg(x'_0,y')$; this follows from
Axiom~\ref{ax:0}. Then, by lifting property for $\Delta\subset
\A^2(K)$, there exists $z'''\in \U$ such that $z'''\simG z''$ and
$(z',z''')\sim_\Delta (x'_0,y')$. Moreover, such a point $z'''$ is
unique. Indeed, by Axiom~\ref{ax:0} the conditions $\phh(z''') \simH
\phh(z'')$ and $(z',z''')\sim_{\Delta,H} (x'_0,y')$ determine
$\phh(z''')$ uniquely for every $H\nrmfin G$. This implies that
$z'''$ is unique by Axiom~\ref{ax:2.1}.

The above construction defines an action $\sigma$ of
$\pi(\U,x'_0)=\{y:y\sim_\Gamma x'_0\}=\pinv p(x'_0)$ on $\U$: a point
$y'\sim_\Gamma x'_0$ sends $z'$ into $z'''$, $\sigma_{y'}z'=z'''$.
Axiom~\ref{ax:0} and Axiom~\ref{ax:1.1} imply that it is in fact a
group action.

Let $\pi(\U)$ be the group of transformations of $\U$ induced by $\pi(U,x'_0)$; the group
does not depend on the choice of $x'_0$. We refer to $\pi(\U)$ as the \emph{group of deck
transformations}, or \emph{the fundamental group} of $\U$. This terminology
is justified by the fact that $\tau\circ p=p$, for $p:\U\ra \A(K)$ the
covering map.

For a subset $W\subset \U^n$, let $\pi(W)=\{\tau:\U ^n\ra \U^n:\tau(W)\subset W, \tau\in\pi(\U)^n\}$.


\subsubsection{Decomposition Lemma for $\U$}

We use a Corollary to Lemma~\ref{lem:noeth:ana}.

\begin{lem}[Decomposition lemma; Noetherian property]\label{lem:noeth:ana:2}\label{lem:decomp.alg}
Assume $\A$ is LERF.

A subset $\pinv(W),W\subset\A(K)$ has a decomposition of the form
$$ W'= H Z'_1\cup \ldots\cup H Z'_k,$$
where $H\nrmfin \Gamma$ is a finite index normal subgroup of $\Gamma$, the
co-etale closed sets $Z'_1,\ldots,Z'_k $ are irreducible components
of relations $\sim_{Z_i}$, for some algebraic subvarieties $Z_i$ of
$\A(K)$, and for any $\tau\in H$ either $\tau Z'_i=Z'_i$ or $\tau
Z'_i\cap Z'_i=\emptyset$.

\end{lem}\bp By a corollary to Decomposition
Lemma~\ref{lem:noeth:ana} we may choose $H\nrmfin\Gamma$
with the following property.

Let $Z_i\subset\etAh(K)$'s be the irreducible components of
$\phh\pinv(W)$. Then, they have the property that the connected
components of $\pgg\phh\inv(Z_i)\subset \etAg(K)$ are irreducible.
Choose $Z'_i$ to be an irreducible components of relations
$\sim_{Z_i}$, i.e.~the closed sets $\phh\inv(Z_i)$.
We claim that these $Z'_i$'s give rise to a decomposition as above.

Before we are able to prove this, let us prove the \emph{lifting
property for $\sim_{Z_i}$}, namely that the map $\phh:Z'_i\ra
Z_i(K)$ is surjective. For convenience, we drop the index $i$ below.

By passing to a smaller subgroup if necessary we may find a variety
$V\subset \etAh(K)^n$ defined over $\barq$ such that for some $g\in
\A^n(K)$, $Z_i$ is a connected component of fibre $V_g$ of $V$ over
$g$, and it holds that if points $x',y'$ are such that $
\phh(x'),\phh(y')\in Z_i$ and $ x'\simW y', \phh(\pr x')=\phh(\pr
y')=g'$ lie in the same connected component of $V'$ over
$g,\phh(g')=g$, then in fact $x'$ and $y'$ lie in the same connected
component of the preimage of $g\times Z_i$, $x'\simZ y'$.

Consider Axiom~\ref{ref:liftfibre} for all $G\nrmfin\Gamma$
sufficiently small \becd
  x' \simc_{V_g,\etAg} y' \lra \exists y'' (y''\sim_G y' \aand
                                    \pr x' = \pr y'' \aand
                                     x' \sim_V y'')
\eecd

Now take any point $ z'\in Z'\subset \U$ and a point $y\in Z(K)$. We
want to prove $\phh(Z')\supset Z(K)$, and thus it is enough to prove
there exists $y_1\in \U$, $\phh(y_1)=y, z'\simZ y_1$. We know that
there exist $y_2\in \U$, $z' \simc_{Z,\etAg} y_2$, due to
Axiom~\ref{ax:0}. Since $Z=V_g$ for some $g\in \U^{n-1}$, we also
have  $(g',z') \simc_{V_g,\etAg} (g',y_2)$, and taking $\phh(g')=g$,
$x'=(g',z'), y'=(g',y_2)$, Axiom~\ref{ref:liftfibre} gives the
conclusion \becd \exists y'' (y'' \sim_G y' \aand
                                    \pr x' = \pr y'' \aand
                                     x' \sim_V y'').\eecd

The conclusion says points $x',y''\in \U^n, \phh(x'),\phh(y')\in Z_i$
lie in the same connected component of $\phh\inv(V)$, are
$\simG$-equivalent, and lie above the same point $g',\phh(g')=g$.
Then by Lemma~\ref{cl:stein} we know that $\phh(x'),\phh(y')$ lie in
the same connected component of the corresponding preimage of $Z_i$.
By definition of $Z'$, this means $\pr_2 y'\in Z'$. Thus, we have
proved that $\phh(Z')=Z(K)$ is surjective.

Now the following by now standard argument concludes the proof.

The the covering property implies that
$$\phh\inv(Z(K))=\bigcup\limits_{h\in H} hZ' = HZ';$$
indeed, by properties of $Z$ we know that the relations
$x'\sim_{Z,G} y'$  are equivalence relations for all $G\nrmfin H$.
Moreover, we know that any two equivalence classes are conjugated by
the action of an element of $H$; this is so because the covering
property implies that there is an element of each of the classes
above each element of $Z(K)$. This implies the lemma. \ep

We single out the following part of the proof as a corollary.

Recall that $\simc$ means ``to lie in the same connected component
of''.
\begin{cor}[the covering property] For a subvariety $Z\subset A(K)$,
 $x' \simc_{Z,G} y' \lra \exists y'' ( y''\simG y' \aand x'\simc_Z
y'').$
\end{cor}
\bp The proof of the lifting property above proves the corollary for
$Z\subset\etAh(K)$ such that the relations $\simc_Z$ and $\simZ$ are
equivalent. However, by Decomposition Lemma any set $\phh\inv(Z)$
can be decomposed into a union of such sets; then going from one
irreducible component to another one intersecting it gives the
corollary.\ep

\begin{cor}[Topology on $\U$]\label{cor:top.gen}
The collection of co-etale closed subsets of $\U$ forms a topology
with a descending chain conditions on irreducible sets. A basic co-etale
closed set possesses an irreducible decomposition as a union of a
finite number of basic co-etale closed sets whose co-etale connected
components are co-etale irreducible. A union of irreducible
components of a co-etale closed set is  co-etale closed.

That is,

\bi \item the collection of co-etale closed subsets on $\uU^n,n>0$
forms a topology. The projection and inclusion maps
$\pr:\uU^n\ra\uU^m, (x_1,\ldots,x_n)\mapsto
(x_{i_1},\ldots,x_{i_m})$ and $\ita:\uU^n\hra\uU^m,
(x_1,\ldots,x_n)\mapsto
(x_{i_1},\ldots,x_{i_{m'}},c_{m'},\ldots,c_{m})$ are continuous.

\item There is no infinite decreasing chain $..\subsetneq
U_{i+1}\subsetneq U_i\subsetneq\ldots\subsetneq U_0$ of 
co-etale closed
irreducible sets.

\item A union of irreducible components of a co-etale closed set is
co-etale closed.

\item A set is basic co-etale closed iff it a union of \emph{connected}
components of a finite number of  $H$-invariant sets, for some
$H\nrmfin\Gamma$ a finite index subgroup of $\Gamma$.

\item A basic co-etale closed set is a union of a finite number of
basic co-etale closed sets whose co-etale connected components are co-etale
irreducible. Moreover, those sets may be taken so that their
connected components within the same set are translates of each
other by the action of a finite index subgroup $H\nrmfin \Gamma$. \ei
\end{cor}
\bp The last item is a reformulation of Decomposition Lemma.
All the items but (1) trivially follow from (5).

Let us prove the intersection of two co-etale closed set $Z'_i$ and
$Y'_i$ is co-etale closed.

Assume $W'$ and $V'$ are unions of connected component of
$H$-invariant sets $H W'$ and $H V'$. The intersection $H W'\cap H
V'$ is $H$-invariant and the set $W'\cap V'$ is a union of the
connected components of $H W'\cap H V'$. The intersection $H W' \cap
HV'=\phh\inv(\phh(W')\cap \phh(V'))$ is co-etale closed by
definition, and thus its connected components are also co-etale
closed. By definition this implies $W'\cap V'$ is co-etale closed.

An infinite intersection is closed by definition. 

The descending chain condition follows from the fact that an
irreducible subset of an irreducible set necessarily has smaller 
dimension.\ep

\subsubsection{Semi-Properness (SP)}

Let $W'\subset \U$ be an \emph{irreducible closed subset} of $\U$,
i.e. a subset of $\U$ defined by $$x\simW a_1 \aand \ldots\aand x\simW
a_n$$ where $a_1,\ldots,a_n\in \U$ are such that $$\forall y \forall
z \left( \AND\limits_{1\leq i\leq n} y \sim a_i \aand
\AND\limits_{1\leq i\leq n} z \sim a_i \lra y\simW z\right).$$ Such
a set $W'$ we call \emph{an irreducible component of closed set
defined by $x\simW x$}, or simply \emph{an irreducible component of
relation $\simW $}.

\begin{lem}[Chevalley Lemma, (SP)] A projection of a co-etale irreducible closed set
is co-etale closed.\end{lem} \bp

Let $W'$ be such an irreducible set, and let $V'=\Cl\pr W'$ be the
least closed set containing its closure. By definition of $V'$
$\phh(\pr W' )\subset \phh(V')$; and by definition of closure
$V'\subset \pr HW'=\phh\inv(\pr \phh(W'))$; the set $\pr\phh(W')$ is
closed by Chevalley Lemma for projective algebraic varieties. The
inequalities imply $\phh(\pr W')=\phh(V')$ for every subgroup
$H\nrmfin\Gamma$.

A deck transformation leaving $W'$ invariant, also leaves $V'$
invariant, i.e.~$\pr\pi(W')\subset\pi(V')$.  On the other hand, the
equality $\phh(\pr W')=\phh(V')$ implies for any $H\nrmfin\Gamma$,
$\pr\pi(W')/H=\pi(V')/H$.

Let us now use Axiom~\ref{ax:2.2} to show that this implies that
$\pr (\pi(W)\cap [H\times H])=\pi(V')\cap H$.

Let us now prove that $\pi(W')\cap H\times H$ is finitely generated
for some $H\nrmfin\Gamma$.

We know by Corollary to Lemma~\ref{cor:gereric.fibres} that $W'=Y'_{g'}$ is a fibre of a $\barq$-defined set $Y'$ over a
point $g'$ such that $\phh(g')\in \pr \phh(Y')$ $\barq$-generic.

We know that for every $G\nrmfin H$, for a connected component $Y_G$
of $\pgg\phh\inv(Y)$, the intersection $Y_G\cap g'\times
\pgg\phh\inv(Y_g)$ is connected; geometrically, that means that a
lifting of $W=Y_g\subset Y$ along the covering map $ Y_G\ra Y$ is a
fibre of $Y$. This holding for every $G\nrmfin H$, it implies that
for $Y'$ a connected component of $\phh\inv(Y)$, the intersection
${Y'}_{g'}^c=Y'\cap g'\times \phh\inv(Y_g)$ is connected, and
therefore it coincides with a connected component of
$\phh\inv(Y_g)=\phh\inv(W)$. Moreover, this implies that if $h\in H$
is such that $h {Y'}_{g'}^c \subset \phh\inv(Y_g)$ then $h {Y'}_{g'}^c
\subset {Y'}_{g'}^c$, i.e.~$h\in\pi({Y'}_{g'}^c)\cap H=\pi(Y'_{g'})\cap
H$. Thus, to prove that $\pi(W)\cap H=\pi({Y'}_{g'}^c)\cap H$ is
finitely generated, it is enough to prove that $\pi(Y'_{g'})\cap H$
is finitely generated. However, the latter is claimed by
Axiom~\ref{ax:2.2} for every variety $Y$ defined over $\barq$.

Let $g_1,\ldots,g_n$ be the generators of $\pi(W')\cap [H\times H]$.
Now take $\tau\in\pi(V')\cap H,\tau(V')=V'$. We know that
$\tau/G\in\pr\pi(W')/G$, for every $G\nrmfin H$, and therefore
$\tau$, up to $\simG$, is expressible as a product of
$g_1,\ldots,g_n$. In other words, that means that $x'$ and $\tau x'$
can be joined by a sequence of points $x'=h_1, h_2,\ldots,h_n=\tau
x'$ such that $h_{i+1}=g_{j_i}h_i$ for all $1\leq i\leq n$, and here
$n=n(G)$ depends on subgroup $G$. By Axiom~\ref{ax:2.3} there is a uniform
bound on such $n=n(G)$, and $\tau$ is expressible as a product of
$g_1,\ldots,g_n$, and therefore belongs to $\pr\pi(W')$.

Now we finish the proof by the covering property argument similar to
the topological proof of Chevalley Lemma in complex case.

  Let $V_0\subset \pr \phh (W') \subset V$ where $V_0\subsetneq V$
  is open in $V$; then $V$ is irreducible. Recall $V'=\Cl\pr W'$ and
  take $V'_0=V'\cap \phh\inv (V_0)$; we know
  $V'_0\subset V'$ is open in $V'$. We also know $V'_0\subset  \Cl\pr W'$.

  Take $v'\in V'_0$, and take $w'\in W'$,
  $\pr\phh(w')=\phh (v')\in V_0\subset \pr W$; such a point $w'$ in $W'$
  exists by the covering property. Now, $\pr w'\in V'$, and thus
  $\gamma_0\in \pi(V')$ where $\gamma_0$ is defined by
  $v'=\gamma_0\pr w'$. Condition $\pr\phh(w')=\phh (v')\in
  \etAh(K)$ implies $\gamma_0\in H$. Thus the inclusion
  $\pr\pi(W')\cap H=\pi(V')\cap H$ implies there exists
  $\gamma_1\in\pi(W')$, $\pr \gamma_1=\gamma_0$, and thus
  $v'=\gamma_0 \pr w' =\pr (\gamma_1w')$, and the Chevalley
  lemma is proven.
\ep



\section{Homogeneity and stability over models}

In the \S\S{} above we have established the main properties of the
co-etale topology on $\U$ (and its Cartesian powers $\U^n$).
That allows us to define and prove the basic properties of
$\Theta$-generic points, for $\Theta$ an algebraically closed
subfield of $K$.

The notion of a $\Theta$-generic point extends to $\U$ in a natural
way. Recall that for a closed $\Theta$-defined set $V'$, the set
$\clt V'$ is the set of all $\Theta$-generic points of $V'$. Recall
also that a set of $\Theta$-generic points of a $\Theta$-defined
set
is called $\Theta$-constructible.

\begin{lem}[Homogeneity]\label{lem:homogeneity}

Any structure $\uu\models\Ax$ is model homogeneous, i.e.~the
projection of a $\Theta$-constructible set is
$\Theta$-constructible, for any algebraically closed subfield
$\Theta$ of the ground field.\end{lem} \bp First note that a point
$w'\in W'$ in an irreducible set $W'$ is $\Theta$-generic iff
$p(w')\in p(W')$ is $\Theta$-generic. By Chevalley Lemma, the fibre
$W'_{g'}$ is non-empty for $g'\in \pr W'$ $\Theta$-generic.
Moreover, by Lemma~\ref{lem:generic.lift} a connected component of
fibre $W_g,g=p(g')$ always contains a $\Theta$-generic point $w\in
W$ of $W$. The lifting $w',p(w')=w$ is always $\Theta$-generic, and
we may find such a lifting in any connected component of a fibre
over a generic point. This implies the lemma.\ep

\begin{defn} Let $\uu,\uu_1,\uu_2\models\Ax$ be $\La$-models of $\Ax(\AC)$
and $\uu\subset\uu_1\cap\uu_2$. We say that tuples $a\in\uu_1^n$ and
$b\in\uu_2^n$ \emph{have the same syntactic quantifier-free type
over $\uu$ in class $\Kk$} if $a$ and $b$ satisfy the same
quantifier-free $\La$-formulae with parameters in $\uu$.
\end{defn}

\begin{defn} A class $\Kk$ of $\La$-structures is
\emph{syntactically stable over countable submodels} iff for any
countable structure $\uu\in \Kk$, the set of complete $\La$-types
over a structure $\uu$ realised in a structure $\uu'\in \Kk$ is at
most countable.\end{defn}

\begin{defn}
 A class $\Kk$ of $\La$-structures is
\emph{quantifier-free syntactically stable over countable submodels} iff there are only countably
many quantifier-free syntactic types in class $\Kk$ over
any countable model $\uu\in\Kk$.
\end{defn}

\begin{lem}[Stability over submodels] Assume $\A$ is LERF.
The class
of $\La$-models of $\Ax(\AC)$ is  quantifier-free syntactically
stable over submodels.\end{lem}
\bp
If $\uu\prec\uu'$ is an elementary substructure, then
$\uu=\uu'(\Theta)=\{u\in\uu': p(u)\in A(\Theta)\}$, for some
algebraically closed subfield $\Theta$.

Every positive quantifier-free $\La$-formula over $\uu$ determines a closed
set defined over $\Theta$.
For every tuple $v'\in \uu'$, there is a least closed set $V'=\Clt(v')$
containing $v'$ and defined over $\Theta$; it is irreducible, and
is a connected component of an algebraic subvariety $V/\Theta$ of
$\etAh$ defined
over $\Theta$, for some $H\nrmfin \Gamma$. Moreover, $\Clt(v')$ has
a $\Theta$-point $v'_\Theta$. Thus, the quantifier-free
$\La$-type of tuple $v'$ is determined by the point
$v'_\Theta\in \uu$ and
a subvariety $V/\Theta$. Therefore, there are only countable number
of such types, which implies that class $\Kk$ is
quantifier-free syntactically stable over submodels.\ep

\begin{thm}[Homogeneity and Stability of class $\Kk$]\label{prp:main}
Assume $\A$ is LERF.

All structures $\La$-models of $\Ax(\AC)$ are model homogeneous.
The class of $\La$-models of $\Ax(\AC)$ is syntactically quantifier-free stable over
countable submodels.\end{thm} \bp Implied by preceding two
lemmata.\ep

Finally, we may state Theorem~\ref{Th:1}, which was the goal of the paper.

\begin{thm}[Model Stability of $\Ax(\uU)$]\label{Th:1} Let $\A$
be a smooth projective algebraic variety which is LERF. 
Let
$\La$ be the countable language defined in Def.~\ref{def:tau}. 
Then

\bi \item[($\aspir{2_{\aleph_0\ra\aleph_1}}$\!)] Any two models
$\uU_1\models\Ax$ and $\uU_2\models \Ax$ of axiomatisation $\Ax$ and
of cardinality $\aleph_1$, such that \bi
    \item[] there exist a common countable submodel
     $\uU_0\models \Ax$, $\uU_0\subset \uU_1$ and $\uU_0\subset\uU_1$
    \ei
are isomorphic, $\uU_1\cong_\La \uU_2$, and, moreover, the
isomorphism $\phi$ is identity on $\uU_0$. \ei
\end{thm}
\bp This is closely related to Proposition~\ref{prp:main}; however,
let us prove this directly in an explicit manner; in this argument
we try to put an emphasis on the properties of the topology,
although this could also be treated as a very common model-theoretic
argument.

We will prove that every partial $\La$-isomorphism $f:\uU_1
\dashrightarrow \uU_2, f(a)=b, a\in \uU_1^n,
f_{|\uU_0}=\id_{|\uU_0}$, $n\in\N$ finite,  defined on $\uU_0\cup
\{a_1,\ldots,a_n\}$, can be extended to $\uU_0\cup
\{a_1,\ldots,a_n\}\cup \{c\}$, $f(c)\in \uU_2$ for any element $c\in
\uU_1$. This allows to extend a partial $\La$-isomorphism from a
{\em countable} model to its {\em countable} extension. This is
enough: by taking unions of chains of countable submodels we get
isomorphism between models of cardinality $\aleph_1$. Note that one
cannot get isomorphism between models of cardinality $\aleph_2$ in
this way.

Let $V_1=\Cl_{\uU_0}(a), W_1=\Cl_{\uU_0}(a,c)$ be the minimal closed
irreducible subsets containing points $a\in\uU_1^n$ and $(a,c)\in\uU_1^{n+1}$; let
$V_2=\Cl_{\uU_0}(f(a))$ be the corresponding subset of $\uU_2$. Since $f$ is an $L$-isomorphism,
sets $V_1$ and $V_2$ are defined by the same $L$-formulae with parameters in $\uU_0$.

Take a subgroup $H\nrmfin \Gamma$ sufficiently small such that $V_1,V_2,W_1,W_2$ are
connected components of $\phh\inv\phh(V_1),\phh\inv\phh(V_2),\phh\inv\phh(W_1),\phh\inv\phh(W_2)$, respectively.
Pick points $v_1,w_1\in\uU_0$ such that $v_1\in V_1,V_2$ and $w_2\in W_1,W_2$.

Now, by definition of $W_2$ we have $\pr \phh W_2=\phh V_2$, and
also $\pr w_2\in V_2$; choose $c'\in \uU_2$ such that
$(\phh(b),\phh(c'))\in \phh(W_2)$ is a $\uU_0$-generic point of
$\phh(W_2)$. Then by the lifting property for $W_2$ there exists a
point $(b',c'')\in W_2$ such that
$\phh(b')=\phh(b),\phh(c'')=\phh(c')$. However, this implies that
$b'\in \pr W_2\subset V_2$ is a $\uU_0$-generic point of $V_2$.
Therefore by the homogeneity properties in
Lemma~\ref{lem:homogeneity}, or equivalently because the projection
$\pr W_2$ is a closed set definable over $\uU_0$, this implies
$V_2\subset \pr W_2$, and, in particular, there exists $d\in \uU_1$
such that $(b,d)\in W_2$ is a $\uU_0$-generic point. Now set
$f(c)=d$. By construction, the points $(a,c)\in \uU_1$ and
$(b,d)\in\uU_2$ lie in the same $\uU_0$-definable closed sets, and,
since every basic relation of $\La$ defines a closed $\uU_0$-defined
set, this implies that $f$ is indeed an $\La$-isomorphism, as
required. \ep

\paragraph*{Acknowledgements}. I thank Boris Zilber and Martin Bays for many useful discussions, 
encouragement and careful proofreading of parts of the manuscript which helped to eliminate 
a few errors and inaccuracies. This paper is closely related to my D Phil thesis, and I wish 
to use this opportunity to express my gratitude to Boris Zilber; the debt to his ideas 
is amply obvious.

%
\bibliographystyle{plain}
\bibliographystyle{abbrv}
\bibliography{nieuw}
\end{document}